\documentclass[12pt]{article}
\usepackage{amsmath}
\usepackage{amsfonts}
\usepackage{amssymb}

\newtheorem{Def}{Definition}[section]
\newtheorem{Prop}{Proposition}[section]
\newtheorem{Le}{Lemma}[section]
\newtheorem{Cor}{Corollary}[section]
\newtheorem{Rem}{Remark}[section]
\newtheorem{Ex}{Example}[section]

\newcommand\ds{\displaystyle}
\newcommand\N{{\mathbb{N}}}
\newcommand\ZR{{\mathbb{Z}}}
\newcommand\R{{\mathbb{R}}}
\newcommand\C{{\mathbb{C}}}

\newcommand\Ac{\mathcal{A}}

\newcommand\Cc{\mathcal{C}}

\newcommand\Ec{\mathcal{E}}
\newcommand\Fc{\mathcal{F}}

\newcommand\Hc{\mathcal{H}}
\newcommand\Kc{\mathcal{K}}
\newcommand\Lc{\mathcal{L}}
\newcommand\Mc{\mathcal{M}}
\newcommand\Nc{\mathcal{N}}
\newcommand\Qc{\mathcal{Q}}
\newcommand\Rc{\mathcal{R}}
\newcommand\Sc{\mathcal{S}}
\newcommand\Tc{\mathcal{T}}

\newcommand\Vc{\mathcal{V}}
\newcommand\Wc{\mathcal{W}}

\newcommand\gc{\mathfrak{g}}

\newcommand\Ker{\mathrm{Ker}}
 
\newcommand\is{i^{*}}
\newcommand\js{j^{*}}
\newcommand\ks{k^{*}}
\newcommand\X[1]{X_{#1}}
\newcommand\Y[1]{Y_{#1}}
\newcommand\Z[1]{Z_{#1}}
\newcommand\W[1]{W_{#1}}
\newcommand\Xh{X^{1,0}}
\newcommand\Yh{Y^{1,0}}
\newcommand\Zh{Z^{1,0}}

\newcommand\Xah{X^{0,1}}
\newcommand\Yah{Y^{0,1}}
\newcommand\Zah{Z^{0,1}}

\newcommand\sg{\sigma}

\newcommand\n[2]{\nu_{#1}^{#2}}

\newcommand\sig[2]{\sg_{#1}^{#2}}

\newcommand\op{\stackrel{\perp}{\oplus}}
\newcommand\tu[1]{\Tc_{#1}}

\newcommand\nab[1]{\nabla_{#1}}

\begin{document}

\title{\LARGE\bf Witt pseudo-Riemannian connection}	   
\author{Robert PETIT*} 
\date{}
\maketitle
\vspace{2cm}
\begin{abstract}
We define the notion of Witt structure on the tangent bundle of a pseudo-Riemannian manifold and we introduce a connection adapted to a such structure. The notions of geodesics and symmetric spaces are revisited in this setting and applications are given in the special cases of Robinson and Fefferman manifolds.
\end{abstract}

\vspace{3cm}\noindent
{\bf Author Keywords:} 
\noindent
{\bf Mathematical subject codes:} 32V05, 53B05, 53B30, 53B35, 53C05, 53C12, 53C17, 53C25, 53C50, 53C55, 53D10, 53Z05, 58E30\\

\vspace{1cm}\noindent
*Laboratoire de Math\'ematiques Jean Leray, UMR 6629 CNRS, Universit\'e de Nantes,\\ 
2, rue de la Houssini\`ere BP 92208, 44322 Nantes - France.\\
E-mail address : robert.petit@univ-nantes.fr

\newpage
\section{Introduction}
Metric connections with torsion play an important part in Hermitian, pseudo-Hermitian and contact Riemannian geometry. Actually, outside the Kahler case, all the Hermitian connections defined on the tangent bundle of an almost Hermitian manifold have torsion. The Lichnerowicz, Chern and Bismut connections are the most famous examples of such connections (cf. \cite{G}). Examples of metric connections with torsion exist in pseudo-Hermitian geometry (Tanaka-Webster connection \cite{NL},\cite{TN},\cite{W}) and in contact Riemannian geometry (Tanno connection \cite{TS}). Recently, in the setting of sub-Riemannian geometry, Diniz-Veloso \cite{DV} and Hladky \cite{H} have introduced for some Riemannian manifolds (obtained as Riemannian extensions of the sub-Riemannian manifolds) endowed with an orthogonal decomposition of the tangent bundle, a metric connection with torsion preserving the decomposition. A substitute to an orthogonal decomposition in the pseudo-Euclidean setting is the well known Witt decomposition. The tangent bundle of some pseudo-Riemannian manifolds like the Robinson manifolds \cite{NT},\cite{TC} or the Fefferman spaces \cite{DT},\cite{F},\cite{L} can be endowed with a type Witt decomposition (called Witt structure in the following). Also, it is possible to define (just as Diniz-Veloso and Hladky) a metric connection with torsion preserving a such decomposition. The main purpose of this article is to construct a such connection and to study its properties. The plan of this article is as follow. In section 2, we give a definition of a Witt structure on the tangent bundle of a pseudo-Riemannian manifold and we exhibit various examples. Section 3 is devoted to construct a metric connection adapted to a Witt structure on a pseudo-Riemannian manifold called the canonical Witt connection (Proposition 3.1). In Section 4, we investigate the notion of geodesic in this new context. In particular, for the Lorentzian manifolds endowed with a Witt structure and a compatible sub-Riemannian structure, we consider the normal sub-Riemannian geodesics (as defined in \cite{S}). The equations of normal sub-Riemannian geodesics are given in terms of the canonical Witt connection (Proposition 4.2 and Corollary 4.1). Note that these equations are similar to those obtained in the sub-Riemannian and pseudo-Hermitian setting (\cite{GGP},\cite{GG},\cite{GGM} and \cite{R}). In Section 5, we introduce a notion of locally Witt symmetric space inspired by the notions of sub-symmetric space (Bieliavski-Falbel-Gorodski \cite{BFG},\cite{FG}) and symmetric Cauchy-Riemann manifold (Kaup-Zaitsev \cite{KZ}). Similar to sub-symmetric case, we show that the notion of Witt symmetric space is related to parallelism properties for torsion and curvature of the canonical Witt connection (Proposition 5.1). Section 6 is dedicated to Robinson manifolds and Fefferman spaces which are Lorentzian conterparts of Hermitian manifolds. By means of the canonical Witt connection associated to the standard complex Witt structure on the complexify tangent bundle of a Robinson manifold, we introduce a connection called Lichnerowicz connection on the real tangent bundle of a Robinson manifold which is the natural conterpart of the Lichnerowicz connection on an almost-Hermitian manifold (Proposition 6.1). The equations of normal sub-Riemannian geodesics (Proposition 6.2) and the parallelism conditions associated to Witt symmetric spaces (Corollary 6.1) are revisited in the special case of Fefferman spaces. Note that we can define various notions of curvature associated to a Witt connection (as the Ricci curvature or the scalar curvature) with the hope to have interesting applications in geometry of pseudo-Riemannian manifolds or in general relativity. This work is still in progress.

\newpage 
\section{Witt pseudo-Riemannian manifolds}
In the following, $(M,g)$ is a pseudo-Riemannian manifold (i.e. a $m$-dimensional manifold $M$ endowed with a non-degenerate symmetric $(0,2)$-tensor field $g$ on $TM$).\\
Let $k,l\in\N^*$, $q_{l}\in\ZR\setminus\N,p_{k},p_{k}^{*}\in\N$. For $r,s\in\N^*$, we set $\ZR_{r}:=\{q_{1},\ldots,q_{r}\}$, $\ZR^{s}=\{p_{1},p_{1}^{*},\ldots,p_{s},p_{s}^{*}\}$ and $\ZR_{r}^{s}:=\ZR_{r}\cup\ZR^{s}$.
\begin{Def}
Let $(M,g)$ be a pseudo-Riemannian manifold.\\A Witt structure on $TM$ is the data of subbundles $\ds{(V_{i})}_{i\in\ZR_{r}^{s}}$ of $TM$ (of rank $m_i$) satisfying \\
i) $\ds TM=(V_{p_{s}}\oplus V_{p_{s}^{*}})\op\ldots\op(V_{p_{1}}\oplus V_{p_{1}^{*}})\op V_{q_{1}}\op\ldots\op V_{q_{r}}$,\\
ii) for any $k\in\{1,\ldots,s\}$, $m_{p_{k}}=m_{p_{k}^{*}}$ and the subbundles $V_{p_{k}}$, $V_{p_{k}^{*}}$ are totally isotropic for $g$ together with $g_{/V_{p_{k}}\times V_{p_{k}^{*}}}$ non-degenerate,\\
iii) for any $l\in\{1,\ldots,r\}$, the subbundles $V_{q_{l}}$ are anisotropic for $g$.
\end{Def}
A Witt pseudo-Riemannian manifold $(M,g,{(V_{i})}_{i\in\ZR_{r}^{s}})$ is a pseudo-Riemannian manifold $(M,g)$ admitting a Witt structure $\ds{(V_{i})}_{i\in\ZR_{r}^{s}}$ on $TM$.

\begin{Rem}
If, for any $k\in\{1,\ldots,s\}$, $V_{p_{k}}\oplus V_{p_{k}^{*}}$ is a hyperbolic plan (i.e. $m_{p_{k}}=m_{p_{k}^{*}}=1$) and $r=1$, then, the Witt structure gives, at each point $x\in M$, a Witt decomposition of $(T_{x}M,g_{x})$.
\end{Rem}

\begin{Def}
A complex Witt structure on a pseudo-Riemannian manifold $(M,g)$ is a Witt structure on the complexification $T^{\C}M$ of $TM$ with respect to the complexified metric $g^{\C}$.
\end{Def}

Now, we give examples of Witt structures on pseudo-Riemannian manifolds.

\begin{Ex} (cf. \cite{BK}, \cite{DOV}, \cite{MR}, \cite{N})
Let $\lambda:=(\lambda_1,\ldots,\lambda_{m})\in\R^{m}$ with $1\leq\lambda_{1}\leq\lambda_{2}\ldots\leq\lambda_{m}$.\\
We respectively denote by $\C_{-}$ and $\C_{+}$ the sets of complex numbers and para-complex numbers (i.e. $\C_{\pm}=\{x+i_{\pm}y,\, x,y\in\R\}$ with $i_{\pm}^{2}=\pm 1$).
We recall that the Heisenberg group $\Hc^{m}$ is the set $\C_{-}^{m}\times\R$ (or $\C_{+}^{m}\times\R$) with the group law 
$$(z_1,z_2,\ldots,z_m,s).(z'_1,z'_2,\ldots,z'_m,s')=(z_1+z'_1,\ldots,z_m+z'_m,s+s'-\frac{1}{2}\sum_{j}Im(z_j\bar{z'_j})).$$
In the following we denote by $\Hc_{-}^{m}$ and $\Hc_{+}^{m}$ depending the choice of $\C_{\pm}$ to define $\Hc^{m}$.\\
Let $\rho_{\pm}^{\lambda}:\R\to Aut(\Hc_{\pm}^{m})$ given by 
$$\rho_{\pm}^{\lambda}(t)(z_1,z_2,\ldots,z_m,s)=(e^{i_{\pm}\lambda_1 t}z_1,\ldots,e^{i_{\pm}\lambda_m t}z_m,s).$$
The oscillator group denoted $osc_{\lambda}$ is the group $osc_{\lambda}:=\Hc_{-}^{m}\rtimes_{\rho_{-}^{\lambda}}\R$.\\
Its Lie algebra $\mathfrak{osc_{\lambda}}=span\{e_{1},\ldots,e_{m},e_{m+1},\ldots,e_{2m},\varepsilon_{0},\varepsilon_{1}\}$ is given by the relations 
$$\begin{array}{ccc}
[\varepsilon_{0},\varepsilon_{1}]=[\varepsilon_{0},e_{i}]=[\varepsilon_{0},e_{m+i}]=0, &&\cr
[\varepsilon_{1},e_{i}]=\lambda_{i}e_{m+i},\quad [\varepsilon_{1},e_{m+i}]=-\lambda_{i}e_{i},&&(i,j\in\{1,\ldots,m\})\cr
[e_{i},e_{j}]=[e_{m+i},e_{m+j}]=0,\quad [e_{i},e_{m+j}]=\delta_{ij}\varepsilon_{0}.&& 
\end{array}$$
We have a Lorentzian bilinear form $\langle\, ,\, \rangle_{\lambda}$ $ad$-invariant on $\mathfrak{osc_{\lambda}}$ given by 
$$\langle x,y\rangle_{\lambda}=\sum_{i=1}^{m}\frac{\lambda_{i}}{2}(x_{i}y_{i}+x_{m+i}y_{m+i})+x_{0}y_{1}+y_{0}x_{1}.$$
With respect to $\langle\, ,\, \rangle_{\lambda}$, we have the decomposition 
$$\mathfrak{osc_{\lambda}}=\mathfrak{(v_{p_{1}}\oplus v_{p_{1}^{*}})\op v_{q_{1}}\op v_{q_{2}}},$$ 
with $\mathfrak{v_{p_{1}}}=span\{\varepsilon_{1}\}$, $\mathfrak{v_{p_{1}^{*}}}=span\{\varepsilon_{0}\}$ and $\mathfrak{v_{q_{1}}}=span\{e_{1},\ldots,e_{m}\}$, $\mathfrak{v_{q_{2}}}=span\{e_{m+1},\ldots,e_{2m}\}$.\\ 
Translating on the left $\langle\, ,\, \rangle_{\lambda}$ and $(e_{1},\ldots,e_{m},e_{m+1},\ldots,e_{2m},\varepsilon_{0},\varepsilon_{1})$, we obtain a bi-invariant Lorentzian metric $g_{\lambda}$ and left-invariant vector fields $(E_{1},\ldots,E_{m},E_{m+1},\ldots,E_{2m},\Ec_{0},\Ec_{1})$ on $osc_{\lambda}$. By setting $V_{p_{1}}=span\{\Ec_{1}\}$, $V_{p_{1}^{*}}=span\{\Ec_{0}\}$, $V_{q_{1}}=span\{E_{1},\ldots,E_{m}\}$, $V_{q_{2}}=span\{E_{m+1},\ldots,E_{2m}\}$, we obtain a Witt structure for the Lorentzian manifold $(osc_{\lambda},g_{\lambda})$  given by 
$$Tosc_{\lambda}=(V_{p_{1}}\oplus V_{p_{1}^{*}})\op V_{q_{1}}\op V_{q_{2}}.$$ 
\\
The hyperbolic oscillator group $osh_{\lambda}$ is the group $osh_{\lambda}:=\Hc_{+}^{m}\rtimes_{\rho_{+}^{\lambda}}\R$.\\
Its Lie algebra $\mathfrak{osh_{\lambda}}=span\{e_{1},\ldots,e_{m},e_{m+1},\ldots,e_{2m},\varepsilon_{0},\varepsilon_{1}\}$ is given by the relations 
$$\begin{array}{ccc}
[\varepsilon_{0},\varepsilon_{1}]=[\varepsilon_{0},e_{i}]=[\varepsilon_{0},e_{m+i}]=0, &&\cr
[\varepsilon_{1},e_{i}]=\lambda_{i}e_{i},\quad [\varepsilon_{1},e_{m+i}]=-\lambda_{i}e_{m+i},&&(i,j\in\{1,\ldots,m\})\cr
[e_{i},e_{j}]=[e_{m+i},e_{m+j}]=0,\quad [e_{i},e_{m+j}]=\delta_{ij}\varepsilon_{0}.&& 
\end{array}$$
A bilinear form of signature $(m+1,m+1)$, $\langle\, ,\, \rangle_{\lambda}$ $ad$-invariant on $\mathfrak{osh_{\lambda}}$ is given by 
$$\langle x,y\rangle_{\lambda}=\sum_{i=1}^{m}\frac{\lambda_{i}}{2}(x_{i}y_{m+i}+x_{m+i}y_{i})+x_{0}y_{1}+y_{0}x_{1}.$$
With respect to $\langle\, ,\, \rangle_{\lambda}$, we have the decomposition 
$$\mathfrak{osh_{\lambda}}=\mathfrak{(v_{p_{2}}\oplus v_{p_{2}^{*}})\op (v_{p_{1}}\oplus v_{p_{1}^{*}})},$$ 
with $\mathfrak{v_{p_{1}}}=span\{\varepsilon_{1}\}$, $\mathfrak{v_{p_{1}^{*}}}=span\{\varepsilon_{0}\}$ and $\mathfrak{v_{p_{2}}}=span\{e_{1},\ldots,e_{m}\}$, $\mathfrak{v_{p_{2}^{*}}}=span\{e_{m+1},\ldots,e_{2m}\}$.\\
Translating on the left, we obtain a bi-invariant pseudo-riemannnian metric $g_{\lambda}$ on $osh_{\lambda}$ and a Witt structure such that 
$$Tosh_{\lambda}=(V_{p_{2}}\oplus V_{p_{2}^{*}})\op (V_{p_{1}}\oplus V_{p_{1}^{*}}).$$
\end{Ex}

\begin{Ex} Almost nul structure (cf. \cite{NT},\cite{P},\cite{TC})\\
An almost nul structure on a $(2m+2+\epsilon)$-dimensional pseudo-Riemannian manifold $(M,g)$ (with $\epsilon\in\{0,1\}$) is a rank-$(m+1)$ complex subbundle $\Nc$ of $T^{\C}M$ which is totally isotropic with respect to the complexified metric $g^{\C}$ (i.e $\Nc\subset\Nc^{\perp}$).\\
In the following, we assume that $\Nc\cap\overline{\Nc}$ (with $\overline{\Nc}$ be the complex conjugate of $\Nc$) is a complex subbundle of $T^{\C}M$ and we denote by $\kappa$ the integer $\kappa:=dim_{\C}(\Nc_{x}\cap\overline{\Nc_{x}})$.\\ 
If $g$ is a riemannian metric then an almost nul structure on $(M,g)$ with $\kappa=0$ is called an almost Hermitian $CR$ structure. An almost Hermitian $CR$ structure provides complex and real Witt structures as follow 
\begin{eqnarray*}
T^{\C}M&=&(\Nc\oplus\overline{\Nc})\op\epsilon{(\Nc\oplus\overline{\Nc})}^{\perp}=(\Vc_{p_{1}}\oplus\Vc_{p_{1}^{*}})\op\Vc_{q_{2}}\\
TM&=&\Qc\op\epsilon{\Qc}^{\perp}=V_{q_{1}}\op V_{q_{2}},
\end{eqnarray*}
where $\Qc$ is the real subbundle of $TM$ given by $\Qc^{\C}=\Nc+\overline{\Nc}$ (with $V^{\C}$ the complexification of $V$).\\
If $g$ is a Lorentzian metric then an almost nul structure on $(M,g)$ with $\kappa=1$ is called an almost Robinson structure. In this case, there exists a nul line subbundle $\Rc$ of $TM$ such that $\Rc^{\C}=\Nc\cap\overline{\Nc}$ and we have the following filtration of $TM$ and $T^{\C}M$ 
$$\{0\}\subset\Rc\subset\Qc\subset\Rc^{\perp}\subset TM,$$
$$\{0\}\subset\Rc^{\C}\subset\Qc^{\C}\subset {\Rc^{\C}}^{\perp}=\Nc^{\perp}+{\overline{\Nc}}^{\perp}\subset T^{\C}M.$$
It follows real and complex Witt structures (associated to the gradings of $TM$ and $T^{\C}M$) 
\begin{eqnarray*}
TM&=&(\Rc\oplus TM/{\Rc}^{\perp})\op(\Qc/\Rc)\op\epsilon(\Rc^{\perp}/\Qc)=(V_{p_{1}}\oplus V_{p_{1}^{*}})\op V_{q_{2}}\op V_{q_{3}},\\
T^{\C}M&=&(\Rc^{\C}\oplus{(TM/{\Rc}^{\perp})}^{\C})\op (\Nc/\Rc^{\C}\oplus\overline{\Nc}/\Rc^{\C})\op\epsilon({(\Rc^{\perp}/\Qc)}^{\C})=(\Vc_{p_{1}}\oplus\Vc_{p_{1}^{*}})\op(\Vc_{p_{2}}\oplus\Vc_{p_{2}^{*}})\op\Vc_{q_{3}}
\end{eqnarray*}
\end{Ex}

\begin{Ex} Fefferman manifolds (cf. \cite{B},\cite{DT},\cite{F},\cite{L})\\
A polycontact form (cf. \cite{VE}) on a smooth manifold $M$ is a $\R^k$-valued $1$-form $\theta$ on $M$ such that, for every $\mu\in{(\R^{k})}^{*}/\{0\}$, the scalar $2$-form $\mu\circ{d\theta}$ is nondegenerate 
on $\Hc_{M}=:\mathrm{Ker}\,\theta$. The subbundle $\Hc_{M}$ of $TM$ is called a polycontact structure. In this case, there exists a unique family $(\xi_{1},\ldots,\xi_{k})$ of vector fields (called Reeb 
fields) satisfying $\theta_{i}(\xi_{j})=\delta_{ij}$ and ${(i(\xi_{i})d\theta_{i})}_{|\Hc_{M}}=0$ and we have the decomposition $TM=\Vc_{M}\oplus\Hc_{M}$ with $\Vc_{M}=span(\xi_{1},\ldots,\xi_{k})$.\\ 
\\
Let $G$ be a Lie group with Lie algebra $\gc$ and $\theta\in{\Omega}^{1}(M;\gc)$ be a $\gc$-valued polycontact form on a manifold $M$. Let $P$ be a principal $G$-bundle over $M$ together with $\pi:P\to M$ the canonical projection and $\sigma\in{\Omega}^{1}(P;\gc)$ be a connection form on $P$. We have the decomposition into horizontal and vertical bundles $TP=\Vc_{P}\oplus\Hc_{P}=\Vc_{P}\oplus\pi^{*}\Vc_{M}\oplus\pi^{*}\Hc_{M}$. Now assume that 
$(M,\Hc_{M},g_{\Hc_{M}})$ is a sub-Riemannian structure and $(\gc,\langle.,.\rangle_{\gc})$ is a metric Lie algebra (Euclidean), then we endow $P$ with the pseudo-Riemannian metric 
$$g^{F}=\pi^{*}g^{\theta}_{\Hc_{M}}+\pi^{*}\theta\odot\sig{}{},$$
where $g^{\theta}_{\Hc_{M}}$ is the extension of $g_{\Hc_{M}}$ to $TM$ such that $\Vc_{M}=\mathrm{Ker}\,g^{\theta}_{\Hc_{M}}$ and where $\odot$ denotes the symmetric product 
(i.e. $\ds(\sigma_{1}\odot\sigma_{2})(X,Y)=\langle\sigma_{1}(X),\sigma_{2}(Y)\rangle_{\gc}+\langle\sigma_{1}(Y),\sigma_{2}(X)\rangle_{\gc}$).\\  
We call $(P,g^{F})$ a Fefferman manifold. Fefferman spaces (defined as $S^1$-principal bundles over strictly pseudoconvex $CR$ manifolds) or contact quaternionic Fefferman spaces (defined as $S^3$-principal bundles 
over contact quaternionic manifolds \cite{A}) are examples of Fefferman manifolds. 
Note that by taking $V_{p_{1}}=\Vc_{P}$, $V_{p_{1}^{*}}=\pi^{*}\Vc_{M}$ and $V_{q_{1}}=\pi^{*}\Hc_{M}$ then $(M,g,{(V_{i})}_{i\in\ZR_{1}^{1}})$ is a Witt pseudo-Riemannian manifold. 
\end{Ex}

\begin{Ex} Lightlike submanifold (cf. \cite{BD}).\\
Let $N$ be a $n$-dimensional submanifold of a pseudo-Riemannian manifold $(M^m,g)$. We assume that, at each point $x\in N$, $rad(T_{x}N))=T_{x}N\cap {T_{x}N}^{\perp}\neq\{0\}$ and that $dim(rad(T_{x}N))$ is constant on $N$ 
($N$ is called lightlike (or degenerate) submanifold of $(M,g)$). Then, the radical bundle $rad(TN)$ is the totally isotropic subbundle of $TM$ whose fibre at $x\in N$ is $rad(T_{x}N)$. Note that $\kappa=rank(rad(TN))$, called the nullity 
degree of $TN$, satisfies $0<\kappa\leq min(n,m-n)$. There exists subbundles $S(TN)\subset TN$, $S(TN^{\perp})\subset TN^{\perp}$ and $ltr(TN)\subset {S(TN^{\perp})}^{\perp}\subset {S(TN)}^{\perp}$ such that 
$$TM=(rad(TN)\oplus ltr(TN))\op S(TN)\op S(TN^{\perp}).$$
The bundles $S(TN)$, $S(TN^{\perp})$ and $ltr(TN)$ are respectively called screen bundle of $TN$, screen bundle of $TN^{\perp}$ and lightlike transversal vector bundle of $N$. By setting 
$V_{p_{1}}=rad(TN)$, $V_{p_{1}^{*}}=ltr(TN)$, $V_{q_{1}}=S(TN)$ and $V_{q_{2}}=S(TN^{\perp})$, then $(M,g,{(V_{i})}_{i\in\ZR_{2}^{1}})$ is a Witt pseudo-Riemannian manifold.
\end{Ex}

\section{Witt connections on Witt pseudo-Riemannian manifolds}
We define in this section a connection adapted to a Witt structure on a pseudo-riemannian manifold. For a connection $\nabla$, the torsion and the curvature of $\nabla$ are respectively defined by  
$$T(X,Y)=\nabla_{X}Y-\nabla_{Y}X-[X,Y]\quad\mathrm{and}\quad R(X,Y)=[\nabla_{X},\nabla_{Y}]-\nabla_{[X,Y]}.$$
In the following, we define $\is$ for $i\in\ZR_{r}^{s}$ by :
$$\begin{array}{ccc}
{\rm if}\, i=q_{l}     &\, {\rm then}\, &\is=q_{l},\cr
{\rm if}\, i=p_{k}     &\, {\rm then}\, &\is=p_{k}^{*},\cr
{\rm if}\, i=p_{k}^{*} &\, {\rm then}\, &\is=p_{k}.
\end{array}
$$
For any point $x\in M$, we denote by $\pi^{i}_{x}:T_{x}M\to V_{i,x}$ the projection onto $V_{i,x}$. Then we have for $X\in T_{x}M$, $X=\ds\sum_{i\in\ZR_{r}^{s}}X_i$ with 
$X_i=\pi^{i}_{x}(X)$. Note that, for any $X,Y\in T_{x}M$, $\ds g(X,Y)=\sum_{i\in\ZR_{r}^{s}}g(X_i,Y_{\is})$.\\

A Witt connection $\nabla$ on a Witt pseudo-Riemannian manifold $(M,g,{(V_{i})}_{i\in\ZR_{r}^{s}})$ is a pseudo-Riemannian connection such that $\nabla(\Gamma(V_{i}))\subset\Gamma(V_{i})$ for any $i\in\ZR_{r}^{s}$.\\
\\
Let $\tu{k}:TM\otimes TM\to V_{k}$ the tensor defined by 
\begin{eqnarray*}
\tu{k}(\X{i},\Y{j})=\left\{ 
\begin{array}{ccc}0 &\ {\rm if}\ & \ i=j \\ 
                  \ds\frac{1}{2}\Bigl({(\Lc_{\X{i}}\Y{j}^{\flat}+\Lc_{\Y{j}}\X{i}^{\flat})}^{\sharp}_{k}-[\X{i},\Y{j}]_{k}-{(d(g(\X{i},\Y{j})))}^{\sharp}_{k}\Bigr)&\ {\rm if}\  &\ i\neq j=k,\ i=\js,\ (j\in\ZR^{s})\\
                  \ds\frac{1}{2}\Bigl({(\Lc_{\X{i}}\Y{j}^{\flat})}^{\sharp}_{k}-[\X{i},\Y{j}]_{k}\Bigr)& \ {\rm if}\ &\ i\neq j=k,\ i\neq\js\\
                  0 &\ {\rm if}\ &\ i\neq j,\ j\neq k
\end{array}\right.
\end{eqnarray*}
with $\X{i}\in\Gamma(V_{i})$, $\Y{j}\in\Gamma(V_{j})$ and where $\flat$,$\sharp$ and $\Lc_{X}$ respectively denote the musical isomorphisms associated to $g$ and the Lie derivative.

\begin{Prop} Let $(M,g,{(V_{i})}_{i\in\ZR_{r}^{s}})$ be a Witt pseudo-Riemannian manifold, then there exists a unique Witt connection $\nabla$ on $TM$ (called the canonical Witt connection) with torsion $T$ given by \\
\begin{equation}\label{e1}
T(\X{i},\Y{j})=-[\X{i},\Y{j}]_{V_{\is\js}^{\perp}}+\tu{j}(\X{i},\Y{j})-\tu{i}(\Y{j},\X{i}),
\end{equation}
where $\ds V_{i_1\ldots i_k}^{\perp}:={(V_{i_1}+\ldots+ V_{i_k})}^{\perp}$.
\end{Prop}

\begin{Rem}
This proposition is the natural extension to the pseudo-Riemannian case of Theorem 3.1 of \cite{DV}and Lemma 2.13 of \cite{H}.
\end{Rem}

We need the following lemma.

\begin{Le} Let $\nabla$ be a pseudo-Riemannian connection, then for any $X,Y,Z\in\Gamma(TM)$, we have :
\begin{equation}\label{e2}
(\nab{X}{Z^{\flat}})(Y)=\frac{1}{2}\Bigl(dZ^{\flat}(X,Y)-Z^{\flat}(T(X,Y))+(\Lc_{Z}g)(X,Y)-g(T(Z,X),Y)-g(T(Z,Y),X)\Bigr),
\end{equation}
\end{Le}
where $Z^{\flat}$ is the $1$-form given by $Z^{\flat}(X)=g(Z,X)$.\\

\noindent{Proof.} For a pseudo-Riemannian connection $\nabla$ and any $X,Y,Z\in\Gamma(TM)$, we have the formula 
\begin{eqnarray}\label{e3}
g(\nab{X}{Y},Z)&=&\frac{1}{2}\Bigl(Xg(Y,Z)-Zg(X,Y)+Yg(X,Z)+g(Z,[X,Y])-g(X,[Y,Z])+g(Y,[Z,X])\Bigr)\nonumber\\
               &&+\Kc(X,Y,Z),
\end{eqnarray}
with $\ds\Kc(X,Y,Z)=\frac{1}{2}\Bigl(g(T(X,Y),Z)-g(T(Y,Z),X)+g(T(Z,X),Y)\Bigr)$ ($\Kc$ is called the contorsion tensor). It follows from (\ref{e3}) that
\begin{eqnarray*}
(\nab{X}{Z^{\flat}})(Y)&=&Xg(Z,Y)-g(Z,\nab{X}{Y})\\
                       &=&\frac{1}{2}\Bigl(Xg(Y,Z)-Yg(X,Z)+Zg(X,Y)-g(Z,[X,Y])+g(X,[Y,Z])\\
                       &&-g(Y,[Z,X])-g(T(X,Y),Z)+g(T(Y,Z),X)-g(T(Z,X),Y)\Bigr)\\
                       &=&\frac{1}{2}\Bigl(XZ^{\flat}(Y)-YZ^{\flat}(X)-Z^{\flat}([X,Y])-Z^{\flat}(T(X,Y))\\
                       &&+Zg(X,Y)-g([Z,X],Y)-g(X,[Z,Y])-g(T(Z,X),Y)-g(T(Z,Y),X)\Bigr).
\end{eqnarray*}
Now, we have   
\begin{eqnarray*}
XZ^{\flat}(Y)-YZ^{\flat}(X)-Z^{\flat}([X,Y])&=&dZ^{\flat}(X,Y)\\
Zg(X,Y)-g([Z,X],Y)-g(X,[Z,Y])&=&(\Lc_{Z}g)(X,Y).
\end{eqnarray*}
The result follows. $\Box$\\

\noindent{Proof of proposition 2.1.} 
For a pseudo-Riemannian connection $\nabla$, the assumption $\nabla$ Witt is equivalent to $(\nab{X}{Z^{\flat}})(\Y{l})=-g(Z,\nab{X}{\Y{l}})=0$ for any $Z\in\Gamma(V_{l}^{\perp})$, $\Y{l}\in\Gamma(V_{l})$ and $X\in\Gamma(TM)$. 
Hence by (\ref{e2}), we have
\begin{equation}\label{e4}
dZ^{\flat}(X,\Y{l})-Z^{\flat}(T(X,\Y{l}))=-(\Lc_{Z}g)(X,\Y{l})+g(T(Z,X),\Y{l})+g(T(Z,\Y{l}),X).
\end{equation}
Let $i,j,k\in\ZR_{r}^{s}$ with $k$ satisfying $k\neq i$ and $k\neq j$, then we have $V_{\ks}\subset V_{i}^{\perp}$ and $V_{\ks}\subset V_{j}^{\perp}$. It follows that $(\nab{\X{i}}{\Z{\ks}^{\flat}})(\Y{j})=(\nab{\Y{j}}{\Z{\ks}^{\flat}})(\X{i})=0$ 
for any $\Z{\ks}\in\Gamma(V_{\ks}),\X{i}\in\Gamma(V_{i}),\Y{j}\in\Gamma(V_{j})$. It follows from these equations together with (\ref{e4}) that \\
for $k\neq i$ and $k\neq j$,
\begin{equation}\label{e5}
\Z{\ks}^{\flat}(T(\X{i},\Y{j})_{k})=d\Z{\ks}^{\flat}(\X{i},\Y{j})\\
\end{equation}
\begin{equation}\label{e6}
g(T(\Z{\ks},\X{i})_{\js},\Y{j})+g(T(\Z{\ks},\Y{j})_{\is},\X{i})=(\Lc_{\Z{\ks}}g)(\X{i},\Y{j}),
\end{equation}
First, since the component $T(\X{i},\Y{i})_{i}$ does not appear in (\ref{e5}) and (\ref{e6}), we take $T(\X{i},\Y{i})_{i}=0$ for any $i\in\ZR_{r}^{s}$.
Now, (\ref{e5}) is equivalent to $g(\Z{\ks},{(T(\X{i},\Y{j})+[\X{i},\Y{j}])}_{k})=0$. By non degeneracy, we deduce that
\begin{equation}\label{e7}
T(\X{i},\Y{j})_{k}=-[\X{i},\Y{j}]_{k} \ {\rm with}\ k\neq i,\ k\neq j.
\end{equation}
If $j\neq\is$ then (\ref{e7}) yields to 
$$
g(T(\Z{\ks},\X{i})_{\js},\Y{j})+g(T(\Z{\ks},\Y{j})_{\is},\X{i})=-g([\Z{\ks},\X{i}]_{\js},\Y{j})-g(\X{i},[\Z{\ks},\Y{j}]_{\is})=(\Lc_{\Z{\ks}}g)(\X{i},\Y{j}).  
$$
Also (\ref{e6}) is always satisfied in this case and is equivalent to 
\begin{eqnarray}\label{e8}
g(T(\Z{\ks},\X{i})_{i},\Y{\is})+g(T(\Z{\ks},\Y{\is})_{\is},\X{i})&=&(\Lc_{\Z{\ks}}g)(\X{i},\Y{\is})\nonumber\\
&=&(\Lc_{\Z{\ks}}\X{i}^{\flat})(\Y{\is})-\Y{\is}^{\flat}([\Z{\ks},\X{i}]_{i}),\quad k\neq i,\ k\neq\is.
\end{eqnarray}
For $i\in\ZR^{s}$, we have $V_{\is}\subset V_{\is}^{\perp}$. Hence, for any $\X{j}\in\Gamma(V_{j})$ with $j\in\ZR_{r}^{s}$, we have $(\nab{\X{j}}{\Z{\is}^{\flat}})(\Y{\is})=0$. 
The case $j\neq i$ is included in equations (\ref{e5}) and (\ref{e6}). Now, if $j=i$, then we have $(\nab{\X{i}}{\Z{\is}^{\flat}})(\Y{\is})=0$ and we obtain, by (\ref{e4}) for $i\in\ZR^{s}$ 
\begin{eqnarray}\label{e9}
g(T(\Z{\is},\X{i})_{i},\Y{\is})-g(T(\Y{\is},\X{i})_{i},\Z{\is})-g(T(\Y{\is},\Z{\is})_{\is},\X{i})&=&g(T(\Z{\is},\X{i})_{i},\Y{\is})-g(T(\Y{\is},\X{i})_{i},\Z{\is})\nonumber\\
                                                                                                 &=&(\Lc_{\Z{\is}}g)(\X{i},\Y{\is})+d\Z{\is}^{\flat}(\X{i},\Y{\is})\nonumber\\
                                                                                                 &=&-(\Lc_{\Y{\is}}g)(\X{i},\Z{\is})-d\Y{\is}^{\flat}(\X{i},\Z{\is})\nonumber\\
                                                                                                 &=&(\Lc_{\Z{\is}}\X{i}^{\flat})(\Y{\is})+(\Lc_{\X{i}}\Z{\is}^{\flat})(\Y{\is})\nonumber\\
                                                                                                 &-&\Y{\is}^{\flat}([\Z{\is},\X{i}]_{i})-(d(g(\Z{\is},\X{i})))(\Y{\is}).
\end{eqnarray}
Now, if we choose, for $k\neq i,\ k\neq\is$,  
\begin{equation}\label{e10}
T(\Z{\ks},\X{i})_{i}=\frac{1}{2}\Bigl({(\Lc_{\Z{\ks}}\X{i}^{\flat})}^{\sharp}_{i}-[\Z{\ks},\X{i}]_{i}\Bigr)
\end{equation}
and, for $i\in\ZR^{s}$,
\begin{equation}\label{e11}
T(\Z{\is},\X{i})_{i}=\frac{1}{2}\Bigl({(\Lc_{\Z{\is}}\X{i}^{\flat}+\Lc_{\X{i}}\Z{\is}^{\flat})}^{\sharp}_{i}-[\Z{\is},\X{i}]_{i}-{(d(g(\Z{\is},\X{i})))}^{\sharp}_{i}\Bigr),
\end{equation}
then (\ref{e8}) and (\ref{e9}) are satisfied. Now formula (\ref{e1}) directly follows from (\ref{e7}),(\ref{e10}),(\ref{e11}) and the assumption $T(\X{i},\Y{i})_{i}=0$. 
By means of the previous formula for $T$ and (\ref{e3}), we define a pseudo-Riemannian connection $\nabla$ with the required properties. $\Box$\\        

\begin{Rem}
Note that the previous proof only requires that $V_{q_{l}}$ are non-degenerate for $l\in\{1,\ldots,r\}$ and not necessarily anisotropic for $g$. 
\end{Rem}

\begin{Ex} Assume that the Witt structure has the following form $(V_{p_{2}}\oplus V_{p_{2}^{*}})\op(V_{p_{1}}\oplus V_{p_{1}^{*}})$ (or $(V_{p_{1}}\oplus V_{p_{1}^{*}})\op V_{q_{1}}\op V_{q_{2}}$) together with $dim\,V_{p_{1}}=1$. If we choose a local field frame $(\n{}{},\n{}{*})$ of $V_{p_{1}}\oplus V_{p_{1}^{*}}$ such that $g(\n{}{},\n{}{*})=1$, then, formula (\ref{e1}) in Proposition 3.1 for the torsion of the canonical Witt connection becomes 
\begin{eqnarray}\label{e12}
T(\X{i},\Y{i})&=&-[\X{i},\Y{i}]_{k}+d\sig{}{}(\X{i},\Y{i})\n{}{}+d\sig{}{*}(\X{i},\Y{i})\n{}{*}\hspace{3.8cm}(k\neq i)\nonumber\\
T(\X{i},\Y{j})&=&d\sig{}{}(\X{i},\Y{j})\n{}{}+d\sig{}{*}(\X{i},\Y{j})\n{}{*}+\tu{j}(\X{i},\Y{j})-\tu{i}(\Y{j},\X{i})\hspace{1.5cm}(i\neq j)\nonumber\\
T(\n{}{},\X{i})&=&-[\n{}{},\X{i}]_{k}+d\sig{}{*}(\n{}{},\X{i})\n{}{*}-\frac{1}{2}(\Lc_{\X{i}}g)(\n{}{},\n{}{*})\n{}{}+\tu{i}(\n{}{},\X{i})\nonumber\\
T(\n{}{*},\X{i})&=&-[\n{}{*},\X{i}]_{k}+d\sig{}{}(\n{}{*},\X{i})\n{}{}-\frac{1}{2}(\Lc_{\X{i}}g)(\n{}{},\n{}{*})\n{}{*}+\tu{i}(\n{}{*},\X{i})\nonumber\\
T(\n{}{},\n{}{*})&=&-[\n{}{},\n{}{*}]_{V_{p_{1}p_{1}^{*}}^{\perp}},
\end{eqnarray}
with $\sig{}{}={\n{}{*}}^{\flat},\;\sig{}{*}={\n{}{}}^{\flat}$ and $i,j,k\in\{p_{2},p_{2}^{*}\}$ (or $i,j,k\in\{q_{1},q_{2}\}$).
In the special case of groups $osc_{\lambda}$ and $osh_{\lambda}$ in example 2.1, then the previous formula becomes 
$$T(\X{i},\Y{i})=0,\quad T(\X{i},\Y{j})=d\sig{}{*}(\X{i},\Y{j})\n{}{*},$$
$$T(\n{}{},Z)=T(\n{}{*},Z)=0.$$
\end{Ex}

\section{Geodesics on Witt Lorentzian manifolds}

In this section, we assume that $(M,g)$ is a Lorentzian manifold endowed with a Witt structure of type $TM=(V_{p_{1}}\oplus V_{p_{1}^{*}})\op V_{q_{1}}$ with $V_{p_{1}},V_{p_{1}^{*}}$ trivial nul line subbundles and $g_{/V_{q_{1}}}$ positive definite on the subbundle $V_{q_{1}}$. We choose (as in section 3) global field frame $(\n{}{},\n{}{*})$ of $V_{p_{1}}\oplus V_{p_{1}^{*}}$ such that $g(\n{}{},\n{}{*})=1$. Denote by $\nabla$ the canonical Witt connection associated to the Witt structure, then the torsion $T$ of $\nabla$ is given by  
\begin{eqnarray}\label{e13}
T(X,Y)&=&d\sig{}{}(X,Y)\n{}{}+d\sig{}{*}(X,Y)\n{}{*},\quad X,Y\in V_{q_{1}}\nonumber\\
T(\n{}{},X)&=&d\sig{}{*}(\n{}{},X)\n{}{*}-\frac{1}{2}(\Lc_{X}g)(\n{}{},\n{}{*})\n{}{}+\tu{V_{q_{1}}}(\n{}{},X)\nonumber\\
T(\n{}{*},X)&=&d\sig{}{}(\n{}{*},X)\n{}{}-\frac{1}{2}(\Lc_{X}g)(\n{}{},\n{}{*})\n{}{*}+\tu{V_{q_{1}}}(\n{}{*},X)\nonumber\\
T(\n{}{},\n{}{*})&=&-[\n{}{},\n{}{*}]_{V_{q_{1}}},
\end{eqnarray}
with $\sig{}{}={\n{}{*}}^{\flat},\;\sig{}{*}={\n{}{}}^{\flat}$ and $\ds g(\tu{V_{q_{1}}}(\n{}{(*)},X),Y)=\frac{1}{2}(\Lc_{\n{}{(*)}}g)(X,Y)$.\\
\\
We define various notions of geodesics adapted to a such type of Witt structure on a Lorentzian manifold.\\ 

{\noindent\bf Lightlike $\nabla$-geodesics.}
\begin{Def} 
A curve $c:[0,1]\to M$ will be called a $\nabla$-geodesic if $\nabla_{\dot{c}}\dot{c}=0$.\\
A lightlike $\nabla$-geodesic is a $\nabla$-geodesic $c:[0,1]\to M$ such that $\dot{c}\in V_{p_{1}}\oplus V_{p_{1}^{*}}$ and $g(\dot{c},\dot{c})=0$.
\end{Def}

\begin{Prop} A curve $c:[0,1]\to M$ is a $\nabla$-lightlike geodesic if and only if
\begin{eqnarray}\label{e14}
(\Lc_{\dot{c}}\sig{}{})(\dot{c})-\sig{}{}(\dot{c})(\Lc_{\n{}{*}}\sig{}{*})(\dot{c})&=&0\nonumber\\
(\Lc_{\dot{c}}\sig{}{*})(\dot{c})-\sig{}{*}(\dot{c})(\Lc_{\n{}{}}\sig{}{})(\dot{c})&=&0.
\end{eqnarray}
\end{Prop}

\noindent{Proof.} Let $c:[0,1]\to M$ with $\dot{c}\in V_{p_{1}}\oplus V_{p_{1}^{*}}$, then $\dot{c}=\sig{}{}(\dot{c})\n{}{}+\sig{}{*}(\dot{c})\n{}{*}$.
Hence, we have 
\begin{eqnarray}\label{e15}
\nabla_{\dot{c}}\dot{c}&=&\nabla_{\dot{c}}(\sig{}{}(\dot{c})\n{}{}+\sig{}{*}(\dot{c})\n{}{*})=(\dot{c}\sig{}{}(\dot{c}))\n{}{}+\sig{}{}(\dot{c})\nabla_{\dot{c}}\n{}{}+(\dot{c}\sig{}{*}(\dot{c}))\n{}{*}+\sig{}{*}(\dot{c})\nabla_{\dot{c}}\n{}{*}\nonumber\\
&=&\bigl((\Lc_{\dot{c}}\sig{}{})(\dot{c})-\sig{}{}(\dot{c})(\nabla_{\dot{c}}\sig{}{})(\n{}{})\bigr)\n{}{}+\bigl((\Lc_{\dot{c}}\sig{}{*})(\dot{c})-\sig{}{*}(\dot{c})(\nabla_{\dot{c}}\sig{}{*})(\n{}{*})\bigr)\n{}{*}.
\end{eqnarray}
Now, using formulas (\ref{e2}) and (\ref{e13}), we have for $X\in TM$ 
$$(\nabla_{X}\sig{}{})(\n{}{})=\frac{1}{2}\Bigl((\Lc_{\n{}{*}}g)(\n{}{},X)-d\sig{}{}(\n{}{},X)-g([\n{}{},\n{}{*}]_{V_{q_{1}}},X)\Bigr).$$
Since $(\Lc_{\n{}{*}}g)(\n{}{},X)=(\Lc_{\n{}{*}}\sig{}{*})(X)+g([\n{}{},\n{}{*}],X)$, we deduce that 
\begin{eqnarray}\label{e16}
(\nabla_{X}\sig{}{})(\n{}{})&=&\frac{1}{2}\Bigl((\Lc_{\n{}{*}}\sig{}{*})(X)-(\Lc_{\n{}{}}\sig{}{})(X)+g([\n{}{},\n{}{*}]-[\n{}{},\n{}{*}]_{V_{q_{1}}},X)\Bigr)\nonumber\\
&=&\frac{1}{2}\Bigl((\Lc_{\n{}{*}}\sig{}{*})(X)-(\Lc_{\n{}{}}\sig{}{})(X)+g(\sig{}{}([\n{}{},\n{}{*}])\n{}{}+\sig{}{*}([\n{}{},\n{}{*}])\n{}{*},X)\Bigr)\nonumber\\
&=&\frac{1}{2}\Bigl((\Lc_{\n{}{*}}\sig{}{*})(X)-(\Lc_{\n{}{}}\sig{}{})(X)+(\Lc_{\n{}{*}}\sig{}{*})(\n{}{})\sig{}{}(X)-(\Lc_{\n{}{}}\sig{}{})(\n{}{*})\sig{}{*}(X)\Bigr)\nonumber\\
&=&\frac{1}{2}\Bigl((\Lc_{\n{}{*}}\sig{}{*})(X+\sig{}{}(X)\n{}{})-(\Lc_{\n{}{}}\sig{}{})(X+\sig{}{*}(X)\n{}{*})\Bigr).
\end{eqnarray}
Applied to $\dot{c}$, (\ref{e16}) yields to 
$$(\nabla_{\dot{c}}\sig{}{})(\n{}{})=(\Lc_{\n{}{*}}\sig{}{*})(\sig{}{}(\dot{c})\n{}{})-(\Lc_{\n{}{}}\sig{}{})(\sig{}{*}(\dot{c})\n{}{*}).$$
Assume that $\dot{c}$ satisfies $g(\dot{c},\dot{c})=0$, then we have $\sig{}{}(\dot{c})\sig{}{*}(\dot{c})=0$ and we obtain    
\begin{equation}\label{e17}
\sig{}{}(\dot{c})(\nabla_{\dot{c}}\sig{}{})(\n{}{})=\sig{}{}(\dot{c})(\Lc_{\n{}{*}}\sig{}{*})(\sig{}{}(\dot{c})\n{}{})=\sig{}{}(\dot{c})(\Lc_{\n{}{*}}\sig{}{*})(\dot{c}).
\end{equation}
In the same way, we have 
\begin{equation}\label{e18}
\sig{}{*}(\dot{c})(\nabla_{\dot{c}}\sig{}{*})(\n{}{*})=\sig{}{*}(\dot{c})(\Lc_{\n{}{}}\sig{}{})(\dot{c}).
\end{equation}
Using (\ref{e17}) and (\ref{e18}), equation (\ref{e15}) becomes 
$$
\nabla_{\dot{c}}\dot{c}=\bigl((\Lc_{\dot{c}}\sig{}{})(\dot{c})-\sig{}{}(\dot{c})(\Lc_{\n{}{*}}\sig{}{*})(\dot{c})\bigr)\n{}{}+\bigl((\Lc_{\dot{c}}\sig{}{*})(\dot{c})-\sig{}{*}(\dot{c})(\Lc_{\n{}{}}\sig{}{})(\dot{c})\bigr)\n{}{*}.
$$
The result directly follows. $\Box$\\
\\
{\bf Normal sub-Riemannian geodesics.}\\
\\
For the Witt structure $TM=(V_{p_{1}}\oplus V_{p_{1}^{*}})\op V_{q_{1}}$, we have $V_{q_{1}}=\mathrm{Ker}\,\sig{}{}\cap\mathrm{Ker}\,\sig{}{*}$ and $g$ split as $g=g_{/V_{q_{1}}}+\sig{}{*}\odot\sig{}{}$. Since $g_{/V_{q_{1}}}$ is a riemannian metric on $V_{q_{1}}$, then by setting $H=V_{q_{1}}$, we obtain a sub-Riemannian structure $(M,H,g_{H}:=g_{/H})$. Recall that a smooth curve $c:[0,1]\to M$ on a sub-riemannian manifold $(M,H,g_{H})$ will be called an horizontal curve if $\dot{c}(t)\in H_{c(t)}$, for any $t\in[0,1]$. In the following, $\Cc$ denotes the set of curves with fixed endpoints (i.e. $\Cc=\{c:[0,1]\to M,\;c(0)=x,c(1)=y\}$) and $\Cc_H$ the subset of horizontal curves with fixed endpoints. Now, we define the length and the energy of a curve respectively by 
$$L_{\Kc}(c)=\int_{0}^{1}\sqrt{g(\dot{c},\dot{c})}\,dt\quad\mathrm{and}\quad E_{\Kc}(c)=\frac{1}{2}\int_{0}^{1}g(\dot{c},\dot{c})\,dt.$$

\begin{Def} 
A normal sub-riemannian geodesic is a critical point of the functional $L_{\Kc}$ on the space $\Cc_H$.
\end{Def}

Up to a reparameterization note that the critical points of $L_{\Kc}$ and $E_{\Kc}$ on $\Cc_H$ coincide. Hence a normal sub-riemannian geodesic is a critical point of the functional $E_{\Kc}$ on the space $\Cc_H$.\\Now, a normal sub-riemannian geodesic can also be viewed as a critical point of a functional $E^{\lambda}$ on $\Cc$ as follow (cf. \cite{GGP}). For $\lambda=(\lambda_{1},\lambda_{2})\in\Cc^{\infty}(M,\R^2)$, we define 
$\alpha^{\lambda}\in\Omega^{1}(M)$ by 
$$\alpha^{\lambda}=\lambda_{1}\sig{}{}+\lambda_{2}\sig{}{*},$$
and $E^{\lambda}:\Cc\to\R$ by 
$$E^{\lambda}(c)=\frac{1}{2}\int_{0}^{1}g(\dot{c},\dot{c})\,dt-\int_{0}^{1}\alpha^{\lambda}(\dot{c})\,dt=E_{\Kc}(c)+\Ac^{\lambda}(c).$$

A curve $c\in\Cc_H$ is a normal sub-riemannian geodesic if and only if there exists $\lambda\in\Cc^{\infty}(M,\R^2)$ such that $c$ is a critical of $E^{\lambda}$ on $\Cc$ ($\lambda$ is the Lagrange multiplier).\\

In order to obtain the variational caracterization of normal sub-riemannian geodesics we recall the concept of variation of a curve. Let $c:[0,1]\to M$ be a smooth curve on a manifold $M$, then a variation of $c$ with fixed endpoints is a smooth map $F:]-\epsilon,\epsilon[\times[0,1]\to M$ such that $F(0,t)=c(t)$ for any $t\in[0,1]$, and, $F(s,0)=c(0)$ and $F(s,1)=c(1)$ for any $s\in]-\epsilon,\epsilon[$.\\ Let $F^{*}TM\to ]-\epsilon,\epsilon[\times[0,1]$ be the pull-back bundle. In the following, we consider the vector fields $w,\dot{c}\in\Gamma(F^{*}TM)$ given by $w_{(s,t)}:=dF_{(s,t)}(\frac{\partial}{\partial s})$ and $\dot{c}_{(s,t)}=dF_{(s,t)}(\frac{\partial}{\partial t})$. Note that $[w_{(s,t)},\dot{c}_{(s,t)}]=0$ and that for a variation with fixed endpoints $w_{(s,0)}=w_{(s,1)}=0$. In the following, a variation of a curve $c$ will be denoted by $c_{s}$ (with $c_{s}(t)=F(s,t)$).

\begin{Prop} Let $c\in\Cc_H$. For any variation $c_{s}$ of $c$, we have
\begin{eqnarray}\label{e19}
\frac{\partial}{\partial s}{E^{\lambda}(c_{s})}_{|s=0}&=&\int_{0}^{1}-g(\nabla_{\dot{c}}\dot{c},w_H)
+\lambda_{1}^{c}d\sig{}{}(\dot{c},w_H)+\lambda_{2}^{c}d\sig{}{*}(\dot{c},w_H)\nonumber\\
&+&\sig{}{}(w)\bigl(\dot{\lambda_{1}^{c}}-\lambda_{1}^{c}(\Lc_{\n{}{}}\sig{}{})(\dot{c})-\lambda_{2}^{c}(\Lc_{\n{}{}}\sig{}{*})(\dot{c})+\frac{1}{2}(\Lc_{\n{}{}}g)(\dot{c},\dot{c})\bigr)\nonumber\\
&+&\sig{}{*}(w)\bigl(\dot{\lambda_{2}^{c}}-\lambda_{1}^{c}(\Lc_{\n{}{*}}\sig{}{})(\dot{c})-\lambda_{2}^{c}(\Lc_{\n{}{*}}\sig{}{*})(\dot{c})+\frac{1}{2}(\Lc_{\n{}{*}}g)(\dot{c},\dot{c})\bigr)\,dt,
\end{eqnarray}
with $\lambda_{i}^{c}:=\lambda_{i}\circ c$ and $w:=w_{(0,t)}$. 
\end{Prop}

\noindent{Proof.} Let $c_{s}$ be a variation of $c$. Then we have
$$E^{\lambda}(c_{s})=E_{\Kc}(c_{s})+\Ac^{\lambda}(c_{s}).$$
Now, 
$$
\frac{\partial}{\partial s}E_{\Kc}(c_{s})=\frac{1}{2}\int_{0}^{1}\frac{\partial}{\partial s}g(\dot{c}_{(s,t)},\dot{c}_{(s,t)})\,dt=\int_{0}^{1}g(\nabla^{F^{*}TM}_{\frac{\partial}{\partial s}}\dot{c}_{(s,t)},\dot{c}_{(s,t)})\,dt,
$$
where $\nabla^{F^{*}TM}$ is the connection on the bundle $F^{*}TM$ induced by the Witt connection $\nabla$ on $TM$.
Now we have the identity 
$$\nabla^{F^{*}TM}_{\frac{\partial}{\partial s}}\dot{c}_{(s,t)}-\nabla^{F^{*}TM}_{\frac{\partial}{\partial t}}w_{(s,t)}-[w_{(s,t)},\dot{c}_{(s,t)}]=T(w_{(s,t)},\dot{c}_{(s,t)}),$$
where $T$ is the torsion of $\nabla$. Since $[w_{(s,t)},\dot{c}_{(s,t)}]=0$, it follows from the previous identity that 
$$
\frac{\partial}{\partial s}E_{\Kc}(c_{s})
=\int_{0}^{1}g(\nabla^{F^{*}TM}_{\frac{\partial}{\partial t}}w_{(s,t)},\dot{c}_{(s,t)})+g(T(w_{(s,t)},\dot{c}_{(s,t)}),\dot{c}_{(s,t)})\,dt.
$$
Using $\ds\frac{\partial}{\partial t}g(w_{(s,t)},\dot{c}_{(s,t)})=g(\nabla^{F^{*}TM}_{\frac{\partial}{\partial t}}w_{(s,t)},\dot{c}_{(s,t)})+g(w_{(s,t)},\nabla^{F^{*}TM}_{\frac{\partial}{\partial t}}\dot{c}_{(s,t)})$ we obtain 
\begin{eqnarray*}
\frac{\partial}{\partial s}E_{\Kc}(c_{s})
&=&\int_{0}^{1}\frac{\partial}{\partial t}g(w_{(s,t)},\dot{c}_{(s,t)})-g(w_{(s,t)},\nabla^{F^{*}TM}_{\frac{\partial}{\partial t}}\dot{c}_{(s,t)})+g(T(w_{(s,t)},\dot{c}_{(s,t)}),\dot{c}_{(s,t)})\,dt\\
&=&g(w_{(s,1)},\dot{c}_{(s,1)})-g(w_{(s,0)},\dot{c}_{(s,0)})+\int_{0}^{1}-g(w_{(s,t)},\nabla^{F^{*}TM}_{\frac{\partial}{\partial t}}\dot{c}_{(s,t)})+g(T(w_{(s,t)},\dot{c}_{(s,t)}),\dot{c}_{(s,t)})\,dt\\
&=&\int_{0}^{1}-g(w_{(s,t)},\nabla^{F^{*}TM}_{\frac{\partial}{\partial t}}\dot{c}_{(s,t)})+g(T(w_{(s,t)},\dot{c}_{(s,t)}),\dot{c}_{(s,t)})\,dt.\\
\end{eqnarray*}
Also we have
$$
\frac{\partial}{\partial s}{E_{\Kc}(c_{s})}_{|s=0}=\int_{0}^{1}-g(w,\nabla_{\dot{c}}\dot{c})+g(T(w,\dot{c}),\dot{c})\,dt.
$$
Using the decomposition $w=w_H+\sig{}{}(w)\n{}{}+\sig{}{*}(w)\n{}{*}$ and formula (\ref{e14}) for the term $g(T(w,\dot{c}),\dot{c})$ with $\dot{c}\in H_{c}$, we obtain 
\begin{equation}\label{e20}
\frac{\partial}{\partial s}{E_{\Kc}(c_{s})}_{|s=0}=\int_{0}^{1}-g(w_{H},\nabla_{\dot{c}}\dot{c})
+\frac{1}{2}\sig{}{}(w)(\Lc_{\n{}{}}g)(\dot{c},\dot{c})+\frac{1}{2}\sig{}{*}(w)(\Lc_{\n{}{*}}g)(\dot{c},\dot{c})\,dt.
\end{equation}
Now, for the second term, we obtain  
$$\frac{\partial}{\partial s}\Ac^{\lambda}(c_{s})=-\int_{0}^{1}\frac{\partial}{\partial s}\bigl(\alpha^{\lambda}(\dot{c}_{(s,t)})\bigr)\,dt
=-\int_{0}^{1}w_{(s,t)}(\alpha^{\lambda}(\dot{c}_{(s,t)}))\,dt.$$
Since $[w_{(s,t)},\dot{c}_{(s,t)}]=0$, we obtain 
$$\frac{\partial}{\partial s}\Ac^{\lambda}(c_{s})=-\int_{0}^{1}(\Lc_{w_{(s,t)}}\alpha^{\lambda})(\dot{c}_{(s,t)})\,dt.$$
By the Cartan formula, we have  
$$(\Lc_{w_{(s,t)}}\alpha^{\lambda})(\dot{c}_{(s,t)})=(d\alpha^{\lambda})(w_{(s,t)},\dot{c}_{(s,t)})+(di(w_{(s,t)})\alpha^{\lambda})(\dot{c}_{(s,t)}).$$
Since $d\alpha^{\lambda}=d\lambda_{1}\wedge\sig{}{}+d\lambda_{2}\wedge\sig{}{*}+\lambda_{1}d\sig{}{}+\lambda_{2}d\sig{}{*}$, we obtain 
\begin{eqnarray*}
(d\alpha^{\lambda})(w_{(s,t)},\dot{c}_{(s,t)})&=&d\lambda_{1}(w_{(s,t)})\sig{}{}(\dot{c}_{(s,t)})-d\lambda_{1}(\dot{c}_{(s,t)})\sig{}{}(w_{(s,t)})+d\lambda_{2}(w_{(s,t)})\sig{}{}(\dot{c}_{(s,t)})-d\lambda_{2}(\dot{c}_{(s,t)})\sig{}{}(w_{(s,t)})\\
&+&\lambda_{1}^{c_{s}}d\sig{}{}(w_{(s,t)},\dot{c}_{(s,t)})+\lambda_{2}^{c_{s}}d\sig{}{*}(w_{(s,t)},\dot{c}_{(s,t)}).
\end{eqnarray*}
and 
$$(di(w_{(s,t)})\alpha^{\lambda})(\dot{c}_{(s,t)})=\frac{\partial}{\partial t}\bigl(i(w_{(s,t)})\alpha^{\lambda}\bigr),$$
with $\lambda_{i}^{c_{s}}=\lambda_{i}\circ c_{s}$. We deduce that 
\begin{eqnarray*}
\frac{\partial}{\partial s}\Ac^{\lambda}(c_{s})&=&-i(w_{(s,1)})\alpha^{\lambda}+i(w_{(s,0)})\alpha^{\lambda}\\
&+&\int_{0}^{1}-d\lambda_{1}(w_{(s,t)})\sig{}{}(\dot{c}_{(s,t)})+d\lambda_{1}(\dot{c}_{(s,t)})\sig{}{}(w_{(s,t)})-d\lambda_{2}(w_{(s,t)})\sig{}{}(\dot{c}_{(s,t)})+d\lambda_{2}(\dot{c}_{(s,t)})\sig{}{}(w_{(s,t)})\\
&-&\lambda_{1}^{c_{s}}d\sig{}{}(w_{(s,t)},\dot{c}_{(s,t)})-\lambda_{2}^{c_{s}}d\sig{}{*}(w_{(s,t)},\dot{c}_{(s,t)}) \,dt\\
&=&\int_{0}^{1}-d\lambda_{1}(w_{(s,t)})\sig{}{}(\dot{c}_{(s,t)})+d\lambda_{1}(\dot{c}_{(s,t)})\sig{}{}(w_{(s,t)})-d\lambda_{2}(w_{(s,t)})\sig{}{}(\dot{c}_{(s,t)})+d\lambda_{2}(\dot{c}_{(s,t)})\sig{}{}(w_{(s,t)})\\
&-&\lambda_{1}^{c_{s}}d\sig{}{}(w_{(s,t)},\dot{c}_{(s,t)})-\lambda_{2}^{c_{s}}d\sig{}{*}(w_{(s,t)},\dot{c}_{(s,t)})\,dt.
\end{eqnarray*}
Using the assumption $\dot{c}(t)=\dot{c}_{(0,t)}\in H$, we have 
\begin{eqnarray*}
\frac{\partial}{\partial s}{\Ac^{\lambda}(c_{s})}_{|s=0}&=&\int_{0}^{1}d\lambda_{1}(\dot{c})\sig{}{}(w)+d\lambda_{2}(\dot{c})\sig{}{*}(w)-\lambda_{1}^{c}d\sig{}{}(w,\dot{c})-\lambda_{2}^{c}d\sig{}{*}(w,\dot{c})\,dt\\
&=&\int_{0}^{1}\dot{\lambda_{1}^{c}}\sig{}{}(w)+\dot{\lambda_{2}^{c}}\sig{}{*}(w)-\lambda_{1}^{c}d\sig{}{}(w,\dot{c})-\lambda_{2}^{c}d\sig{}{*}(w,\dot{c})\,dt.
\end{eqnarray*}
Since $w=w_H+\sig{}{}(w)\n{}{}+\sig{}{*}(w)\n{}{*}$, we deduce that  
\begin{eqnarray}\label{e21}
\frac{\partial}{\partial s}{\Ac^{\lambda}(c_{s})}_{|s=0}&=&\int_{0}^{1}-\lambda_{1}^{c}d\sig{}{}(w_H,\dot{c})-\lambda_{2}^{c}d\sig{}{*}(w_H,\dot{c})+\sig{}{}(w)\bigl(\dot{\lambda_{1}^{c}}-\lambda_{1}^{c}d\sig{}{}(\n{}{},\dot{c})-\lambda_{2}^{c}d\sig{}{*}(\n{}{},\dot{c})\bigr)\nonumber\\
&&+\sig{}{*}(w)\bigl(\dot{\lambda_{2}^{c}}-\lambda_{1}^{c}d\sig{}{}(\n{}{*},\dot{c})-\lambda_{2}^{c}d\sig{}{*}(\n{}{*},\dot{c})\bigr)\,dt\nonumber\\
&=&\int_{0}^{1}-\lambda_{1}^{c}d\sig{}{}(w_H,\dot{c})-\lambda_{2}^{c}d\sig{}{*}(w_H,\dot{c})+\sig{}{}(w)\bigl(\dot{\lambda_{1}^{c}}-\lambda_{1}^{c}(\Lc_{\n{}{}}\sig{}{})(\dot{c})-\lambda_{2}^{c}(\Lc_{\n{}{}}\sig{}{*})(\dot{c})\bigr)\nonumber\\
&+&\sig{}{*}(w)\bigl(\dot{\lambda_{2}^{c}}-\lambda_{1}^{c}(\Lc_{\n{}{*}}\sig{}{})(\dot{c})-\lambda_{2}^{c}(\Lc_{\n{}{*}}\sig{}{*})(\dot{c})\bigr)\,dt.
\end{eqnarray}
The result follows from (\ref{e20}) and (\ref{e21}). $\Box$\\

We deduce from (\ref{e19}) that  

\begin{Cor} A curve $c\in\Cc_H$ is a normal sub-riemannian geodesic if and only if 
\begin{eqnarray}\label{e22}
\nabla_{\dot{c}}\dot{c}&=&\lambda_{1}^{c}{\bigl(i(\dot{c})d\sig{}{}\bigr)}_{H}^{\sharp}+\lambda_{2}^{c}{\bigl(i(\dot{c})d\sig{}{*}\bigr)}_{H}^{\sharp}\nonumber\\
\begin{pmatrix}\dot{\lambda_{1}^{c}}\cr\dot{\lambda_{2}^{c}}\end{pmatrix}&=&\begin{pmatrix}(\Lc_{\n{}{}}\sig{}{})(\dot{c}) & (\Lc_{\n{}{}}\sig{}{*})(\dot{c})\cr
(\Lc_{\n{}{*}}\sig{}{})(\dot{c}) & (\Lc_{\n{}{*}}\sig{}{*})(\dot{c})
\end{pmatrix}\begin{pmatrix}\lambda_{1}^{c}\cr\lambda_{2}^{c}\end{pmatrix}+\frac{1}{2}\begin{pmatrix}(\Lc_{\n{}{}}g)(\dot{c},\dot{c})\cr(\Lc_{\n{}{*}}g)(\dot{c},\dot{c})\end{pmatrix}.
\end{eqnarray}
\end{Cor}

\section{Curvature of Witt connections and Witt symmetric spaces}

\noindent{\bf Curvature of Witt connections}\\

Let $(M,g,{(V_{i})}_{i\in\ZR_{r}^{s}})$ be a Witt pseudo-Riemannian manifold endowed with its canonical Witt connection $\nabla$. Let $R\in{\Omega}^{2}(M;End(TM))$ be the curvature of $\nabla$. Then recall that $R$ satisfies the Bianchi identities 
\begin{eqnarray*}
b(R)(X,Y,Z)&=&(d^{\nabla}T)(X,Y,Z)\\
(d^{\nabla}R)(X,Y,Z)&=&0.
\end{eqnarray*}
where $d^{\nabla}$ is the exterior covariant derivative on vector valued $2$-forms on $M$ and $b(R)$ is the Bianchi map (i.e. $b(R)(X,Y,Z)=R(X,Y)Z+R(Z,X)Y+R(Y,Z)X$). 
The curvature tensor of the canonical Witt connection is given by 
$$R(X,Y,Z,W)=g(R(X,Y)Z,W).$$
Since $\nabla$ is metric, $R$ satisfies the identities $R(Y,X,Z,W)=-R(X,Y,Z,W)$ and $R(X,Y,W,Z)=-R(X,Y,Z,W)$ but $R$ is not symmetric and not satisfies the first Bianchi identity (since $\nabla$ has torsion). 
It directly follows from the first Bianchi identity that 
\begin{eqnarray*}
R(X,Y)\Z{i}&=&{((d^{\nabla}T)(X,Y,\Z{i}))}_{i},\hspace{1.5cm}X,Y\in\Gamma(V_{\is}^{\perp}),\hspace{0.5cm}\X{i},\Z{i}\in\Gamma(V_{i}),\\
R(\X{i},Y)\Z{i}-R(\Z{i},Y)\X{i}&=&{((d^{\nabla}T)(\X{i},Y,\Z{i}))}_{i}.
\end{eqnarray*}

\noindent{\bf Witt symmetric spaces}\\

Following \cite{BFG} and \cite{KZ}, we now introduce a notion of symmetric space adapted to a Witt structure on a pseudo-Riemannian manifold.\\

Let $M$ be a manifold endowed with an affine connection $\nabla$ and $\psi:M\to M$ be a diffeomorphism.
We defined the connection $\psi^{*}\nabla$ on $M$ by setting 
$$(\psi^{*}\nabla)_{X}Y:=\psi^{*}(\nabla_{\psi_{*}(X)}\psi_{*}(Y))=d\psi^{-1}(\nabla_{d\psi(X)}d\psi(Y)).$$
We recall that $\psi$ is an affine diffeomorphism if $(\psi^{*}\nabla)_{X}Y=\nabla_{X}Y$ (equivalently $(\nab{X}d\psi)(Y)=\nab{d\psi(X)}d\psi(Y)-d\psi(\nab{X}Y)=0$).
Note that any affine diffeomorphism preserves torsion and curvature of the connection $\nabla$ (and their derivatives). 

\begin{Def} Let $(M,g,{(V_{i})}_{i\in\ZR_{r}^{s}})$ be a Witt pseudo-Riemannian manifold. A map $\psi:M\to M$ is said to be Witt at a point $x\in M$ if, $d\psi(x)\circ\pi^{i}_{x}=\pi^{i}_{\psi(x)}\circ d\psi(x)$ for any $i\in\ZR_{r}^{s}$ (i.e. $d\psi(x)(V_{i,x})\subset V_{i,\psi(x)}$ for any $i\in\ZR_{r}^{s}$). 
A map $\psi:M\to M$ is said to be Witt if it is Witt at each point $x\in M$. 
\end{Def}

\begin{Def} Let $(M,g,{(V_{i})}_{i\in\ZR_{r}^{s}})$ be a Witt pseudo-Riemannian manifold endowed with its canonical Witt connection. A (local) Witt affine diffeomorphism on $M$ is a (local) Witt affine diffeomorphism on $M$.
A (local) Witt isometry on $M$ is a (local) Witt isometric diffeomorphism on $M$. 
\end{Def}

\begin{Le} Let $(M,g,{(V_{i})}_{i\in\ZR_{r}^{s}},\nabla)$ be a Witt pseudo-Riemannian manifold endowed with its canonical Witt connection, then any local Witt isometry is a local Witt affine diffeomorphism.
\end{Le}

\noindent{Proof.} Let $\psi$ be a local Witt isometry. We have to proove that $g((\nab{X}d\psi)(Y),d\psi(Z))=0$ for any $X,Y,Z\in\Gamma(TM)$. Since $d\psi(V_j)\subset V_j$ and $\nabla$ preserves $V_j$, then we obtain  
$$g((\nab{X}d\psi)(Y),d\psi(Z))=\sum_{j\in\ZR_{r}^{s}}g((\nab{X}d\psi)(\Y{j}),d\psi(\Z{\js})),$$
with $\Y{j}\in\Gamma(V_j),\Z{\js}\in\Gamma(V_{\js})$. Hence, it is sufficient to verify that $\ds g((\nab{X}d\psi)(\Y{j}),d\psi(\Z{\js}))=0$ for $j\in\ZR_{r}^{s}$. 
For any symmetric $2$-tensor $h$ on $M$ and any map $\psi:M\to M$ such that $\nabla h=\nabla\psi^{*}h=0$, we have the following formula 
\begin{eqnarray*}
h((\nab{X}d\psi)(Y),d\psi(Z))&=&\frac{1}{2}\Bigl(h(T(d\psi(X),d\psi(Y)),d\psi(Z))-h(T(d\psi(X),d\psi(Z)),d\psi(Y))\\
                             &&-h(T(d\psi(Y),d\psi(Z)),d\psi(X))\\
                             &&-\psi^{*}h(T(X,Y),Z)+\psi^{*}h(T(X,Z),Y)+\psi^{*}h(T(Y,Z),X)\Bigr).
\end{eqnarray*}
If $\psi$ is a Witt isometry, then the previous formula applied with $h=g$ yields for $\X{i}\in\Gamma(V_{i})$ 
\begin{eqnarray}\label{e23}
g((\nab{\X{i}}d\psi)(\Y{j}),d\psi(\Z{\js}))&=&\frac{1}{2}\Bigl(g(T(d\psi(\X{i}),d\psi(\Y{j}))_{j},d\psi(\Z{\js}))-g(T(d\psi(\X{i}),d\psi(\Z{\js}))_{\js},d\psi(\Y{j}))\nonumber\\
                                   &&-g(T(d\psi(\Y{j}),d\psi(\Z{\js}))_{\is},d\psi(\X{i}))\nonumber\\
                                   &&-g(T(\X{i},\Y{j})_{j},\Z{\js})+g(T(\X{i},\Z{\js})_{\js},\Y{j})+g(T(\Y{j},\Z{\js})_{\is},\X{i})\Bigr).
\end{eqnarray}
If $i=j$ (respectively $i=\js$), then (\ref{e23}) together with the assumption $T(\X{j},\Y{j})_{j}=0$, yields 
\begin{eqnarray}\label{e24}
g((\nab{\X{j}}d\psi)(\Y{j}),d\psi(\Z{\js}))&=&\frac{1}{2}\Bigl(-g(T(d\psi(\X{j}),d\psi(\Z{\js}))_{\js},d\psi(\Y{j}))\nonumber\\
                                           &&-g(T(d\psi(\Y{j}),d\psi(\Z{\js}))_{\js},d\psi(\X{j}))\nonumber\\
                                           &&+g(T(\X{j},\Z{\js})_{\js},\Y{j})+g(T(\Y{j},\Z{\js})_{\js},\X{j})\Bigr),
\end{eqnarray}
(respectively
\begin{eqnarray}\label{e25}
g((\nab{\X{\js}}d\psi)(\Y{j}),d\psi(\Z{\js}))&=&\frac{1}{2}\Bigl(g(T(d\psi(\X{\js}),d\psi(\Y{j}))_{j},d\psi(\Z{\js}))\nonumber\\
                                             &&+g(T(d\psi(\Z{\js}),d\psi(\Y{j}))_{j},d\psi(\X{\js}))\nonumber\\
                                             &&-g(T(\X{\js},\Y{j})_{j},\Z{\js})-g(T(\Z{\js},\Y{j})_{j},\X{\js})\Bigr)).
\end{eqnarray}
If $j\in\ZR_{r}$, then $j=\js$ and $g((\nab{\X{j}}d\psi)(\Y{j}),d\psi(\Z{j}))=0$ since $T(\X{j},\Y{j})_{j}=0$.\\ 
If $j\in\ZR^{s}$, we have, using properties of $T$, that $\ds g(T(\X{j},\Z{\js})_{\js},\Y{j})+g(T(\Y{j},\Z{\js})_{\js},\X{j})=0$. Also, we obtain by (\ref{e24}),(\ref{e25}) that 
$g((\nab{\X{j}}d\psi)(\Y{j}),d\psi(\Z{\js}))=g((\nab{\X{\js}}d\psi)(\Y{j}),d\psi(\Z{\js}))=0$.\\
Now, if $i\neq j$ and $i\neq\js$ then, using properties of $T$, we obtain that 
$$g(T(d\psi(\X{i}),d\psi(\Y{j}))_{j},d\psi(\Z{\js}))-g(T(d\psi(\X{i}),d\psi(\Z{\js}))_{\js},d\psi(\Y{j}))=0$$
and 
$$g(T(\X{i},\Y{j})_{j},\Z{\js})-g(T(\X{i},\Z{\js})_{\js},\Y{j})=0.$$
Hence, we obtain by (\ref{e23})
\begin{eqnarray*}
g((\nab{\X{i}}d\psi)(\Y{j}),d\psi(\Z{\js}))&=&\frac{1}{2}\Bigl(-{(d\psi(\X{i}))}^{\flat}(T(d\psi(\Y{j}),d\psi(\Z{\js}))_{\is})+{\X{i}}^{\flat}(T(\Y{j},\Z{\js})_{\is})\Bigr)\\
&=&-\frac{1}{2}\Bigl(d{(d\psi(\X{i}))}^{\flat}(d\psi(\Y{j}),d\psi(\Z{\js}))-d{\X{i}}^{\flat}(\Y{j},\Z{\js})\Bigr)\\
&=&-\frac{1}{2}\Bigl(({\psi}^{*}d{(d\psi(\X{i}))}^{\flat})(\Y{j},\Z{\js})-d{\X{i}}^{\flat}(\Y{j},\Z{\js})\Bigr)\\
&=&-\frac{1}{2}\Bigl(d({\psi}^{*}{(d\psi(\X{i}))}^{\flat})(\Y{j},\Z{\js})-d{\X{i}}^{\flat}(\Y{j},\Z{\js})\Bigr).
\end{eqnarray*}
Since $\psi$ is a local isometry then ${\psi}^{*}{(d\psi(\X{i}))}^{\flat}={\X{i}}^{\flat}$ and $g((\nab{\X{i}}d\psi)(\Y{j}),d\psi(\Z{\js}))=0$. 
Hence $\psi$ is a local affine diffeomorphism. $\Box$\\

\begin{Rem}
It follows that any Witt isometry maps $\nabla$-geodesics to $\nabla$-geodesics. Also if $\psi_{1},\psi_{2}$ are two Witt isometries on a connected Witt pseudo-Riemannian manifold such that 
$\psi_{1}(x)=\psi_{2}(x)$ and $d\psi_{1}(x)=d\psi_{2}(x)$ then $\psi_{1}=\psi_{2}$.
\end{Rem}

From now, we assume that $\ZR_{r}^{s}$ contains even and odd integers, and, that for any $k\in\{1,\ldots,s\}$, $p_{k}=p_{k}^{*}\;{\rm mod}\,2$. It follows that the involution $\delta_{x}:T_{x}M\to T_{x}M$ 
given by $\ds\delta_{x}=\sum_{i\in\ZR_{r}^{s}}{(-1)}^{i}\pi^{i}_{x}\;(\neq id_{x})$ is a Witt linear isometry. Let $exp^{\nabla}:TM\to M$ be the exponential map associated to $\nabla$ 
then, for any $x\in M$, there exists $U_{x}\subset M$ and $W_{0_{x}}\subset T_{x}M$ neighbourhoods of $x\in M$ and $0_{x}\in T_{x}M$ such that $exp^{\nabla}_{x}:W_{0_{x}}\to U_{x}$ be a diffeomorphism. 
Hence, we can define a Witt diffeomorphism $\psi_{x}:U_{x}\to U_{x}$ called local Witt symmetry at $x$ by 
$$\psi_{x}=exp^{\nabla}_{x}\circ\delta_{x}\circ {(exp^{\nabla}_{x})}^{-1}.$$ 
Note that $\psi_{x}^{2}=id_{M}$, $\psi_{x}(x)=x$ and ${\psi_{x}}_{*}(x)=d\psi_{x}(x)=\delta_{x}$.

\begin{Def} A (locally) Witt symmetric space is a Witt pseudo-Riemannian manifold $(M,g,{(V_{i})}_{i\in\ZR_{r}^{s}},\nabla)$ such that, for any point $x\in M$, $\psi_{x}$ is a (local) Witt isometry.
\end{Def}

\begin{Rem}
It follows from properties of $\psi_{x}$ and remark 4.1. that $\psi_{x}$ is unique.
\end{Rem}

In the following, $\ds V_{\pm}\to M$ are the vector bundles over $M$ with fibre 
$$\ds V_{\pm,x}=ker(\delta_{x}\mp id_{x})=\bigoplus_{i\; even/odd}V_{i,x}.$$

\begin{Prop} Let $(M,g,{(V_{i})}_{i\in\ZR_{r}^{s}},\nabla)$ be a locally Witt symmetric space, then we have
$$[V_{+},V_{+}]\subset V_{+},\quad T_{/V_{-}\times V_{-}}\subset V_{+},\quad T_{/V_{+}\times V_{-}}\subset V_{-},\quad R_{/V_{+}\times V_{-}}=0,$$
and
$$\nab{V_{-}}R=\nab{V_{-}}T=0.$$
\end{Prop}

\noindent{Proof.}
Let $x_{0}\in M$ and $\psi$ be the local Witt symmetry at $x_{0}$. Since $\psi$ is a Witt affine diffeomorphism, we deduce that for $\X{i}\in V_{i,x_{0}},\ \Y{j}\in V_{j,x_{0}},\ \Z{k}\in V_{k,x_{0}}$
$$T(d\psi(\X{i}),d\psi(\Y{j}))=d\psi(T(\X{i},\Y{j}))$$
and
$$R(d\psi(\X{i}),d\psi(\Y{j}))d\psi(\Z{k})=d\psi(R(\X{i},\Y{j})\Z{k}).$$
If $i+j$ is even (respectively odd) we obtain 
\begin{equation}\label{e26}
\sum_{l\in\ZR_{r}^{s}}{(T(\X{i},\Y{j}))}_{l}=\sum_{l\in\ZR_{r}^{s}}{(-1)}^{l}{(T(\X{i},\Y{j}))}_{l},
\end{equation}
(respectively
\begin{eqnarray}\label{e27}
-\sum_{l\in\ZR_{r}^{s}}{(T(\X{i},\Y{j}))}_{l}=\sum_{l\in\ZR_{r}^{s}}{(-1)}^{l}{(T(\X{i},\Y{j}))}_{l},\nonumber\\
-R(\X{i},\Y{j})\Z{k}=R(\X{i},\Y{j})\Z{k}).
\end{eqnarray}
It directly follows from (\ref{e26}) that if $i+j$ is even then $\ds{(T(\X{i},\Y{j}))}_{2l+1}=0$. In addition, if $i,j$ are both even, then $\ds{(T(\X{i},\Y{j}))}_{2l+1}=-{[\X{i},\Y{j}]}_{2l+1}=0$ and consequently
$[V_{+},V_{+}]\subset V_{+}$. If $i+j$ is odd then (\ref{e27}) yields to $\ds{(T(\X{i},\Y{j}))}_{2l}=0$ and $R(\X{i},\Y{j})=0$.
Now, the assumption $\psi$ Witt affine diffeomorphism also implies that for $\X{k}\in V_{k,x_{0}},\ \Y{i}\in V_{i,x_{0}},\ \Z{j}\in V_{j,x_{0}}\ {\rm and}\ \W{l}\in V_{l,x_{0}}$ then we have 
$$(\nab{d\psi(\X{k})}T)(d\psi(\Y{i}),d\psi(\Z{j}))=d\psi((\nab{\X{k}}T)(\Y{i},\Z{j})),$$
and
$$(\nab{d\psi(\X{k})}R)(d\psi(\Y{i}),d\psi(\Z{j}))d\psi(\W{l})=d\psi((\nab{\X{k}}R)(\Y{i},\Z{j})\W{l}).$$
If $i+j$ is even (respectively $i+j$ is odd) and $k$ is odd then we obtain 
\begin{eqnarray}\label{e28}
-\sum_{l\in\ZR_{r}^{s}}{((\nab{\X{k}}T)(\Y{i},\Z{j}))}_{l}&=&\sum_{l\in\ZR_{r}^{s}}{(-1)}^{l}{((\nab{\X{k}}T)(\Y{i},\Z{j}))}_{l},\nonumber\\
-(\nab{\X{k}}R)(\Y{i},\Z{j})\W{l}&=&(\nab{\X{k}}R)(\Y{i},\Z{j})\W{l},
\end{eqnarray}
(respectively
\begin{equation}\label{e29}
\sum_{l\in\ZR_{r}^{s}}{((\nab{\X{k}}T)(\Y{i},\Z{j}))}_{l}=\sum_{l\in\ZR_{r}^{s}}{(-1)}^{l}{((\nab{\X{k}}T)(\Y{i},\Z{j}))}_{l})\;.
\end{equation}
It follows from (\ref{e28}) and (\ref{e29}) that if $i+j$ is even (respectively odd) and $k$ is odd then $\ds(\nab{\X{k}}R)(\Y{i},\Z{j})=0$ and $\ds{((\nab{\X{k}}T)(\Y{i},\Z{j}))}_{2l}=0$ 
(respectively $\ds{((\nab{\X{k}}T)(\Y{i},\Z{j}))}_{2l+1}=0$). We deduce that if $i+j$ is even and $k$ is odd then $\ds(\nab{\X{k}}T)(\Y{i},\Z{j})=0$ and $\ds(\nab{\X{k}}R)(\Y{i},\Z{j})=0$. 
If $i+j$ is odd then $R(\Y{i},\Z{j})=0$ and, together with $k$ odd, then $\ds(\nab{\X{k}}T)(\Y{i},\Z{j})=0$. $\Box$

\begin{Rem}
A bracket generating condition on the distribution $V_{-}$ seems to be necessary to obtain a converse of the proposition. 
\end{Rem}

\section{Robinson manifolds and Fefferman spaces}
Let $(M^{2m+2},g,\Nc)$ be an even dimensional Lorentzian manifold endowed with an almost Robinson structure (since dim $M$ even then  $\Nc^{\perp}=\Nc$). We call Robinson structure an almost Robinson structure such that $\Nc$ is integrable (i.e. $[\Gamma(\Nc),\Gamma(\Nc)]\subset\Gamma(\Nc)$). In the following, we call almost Robinson manifold (respectively Robinson manifold) an even dimensional Lorentzian manifold endowed with an almost Robinson structure (respectively Robinson structure).\\
Denote by $\Rc$ the nul line subbundle of $TM$ such that $\Rc^{\C}=\Nc\cap\overline{\Nc}$. From now, we assume that $(M^{2m+2},g,\Nc)$ is an almost Robinson manifold together with $\Rc$ trivial. Fix a nul line subbundle ${\Rc}^{*}$ of $TM$ (dual to $\Rc$) such that $TM={\Rc}^{*}\oplus{\Rc}^{\perp}$ and consider $\n{}{}$ and $\n{}{*}$ nowwhere vanishing sections of $\Rc$ and ${\Rc}^{*}$ such that $g(\n{}{},\n{}{*})=1$. We have 
${\Rc}^{\perp}=\mathrm{Ker}\,\sig{}{*}$ and ${{\Rc}^{*}}^{\perp}=\mathrm{Ker}\,\sig{}{}$ with $\sig{}{}={\n{}{*}}^{\flat}$ and $\sig{}{*}={\n{}{}}^{\flat}$. Let $\Sc$ the subbundle of $TM$ (called screen bundle) given by $\Sc={(\Rc\oplus{\Rc}^{*})}^{\perp}(=\mathrm{Ker}\,\sig{}{}\cap\mathrm{Ker}\,\sig{}{*})$. Note that ${\Rc}^{\perp}=\Rc\oplus\Sc$ and (as in section 4) that $g_{/\Sc}$ is a riemannian metric on $\Sc$.
Now, there exists a complex subbundle of $T^{\C}M$ denoted by $T^{1,0}M$ such that 
$$\Nc=\Rc^{\C}\oplus T^{1,0}M\;\mathrm{and}\;\Sc^{\C}=T^{1,0}M\oplus T^{0,1}M,$$ 
with $T^{0,1}M=\overline{T^{1,0}M}$. Hence, we have the complex and real Witt structures 
\begin{eqnarray*}
T^{\C}M&=&(\Rc^{\C}\oplus{(\Rc^{*})}^{\C})\op(T^{1,0}M\oplus T^{0,1}M)\\
TM&=&(\Rc\oplus{\Rc}^{*})\op\Sc.
\end{eqnarray*}

We can define an almost Hermitian structure $(g_{\Sc},J)$ on $\Sc$ such that $T^{1,0}M$ (respectively $T^{0,1}M$) is the $\sqrt{-1}$ (respectively $-\sqrt{-1}$)-eigenbundle of $J$ on $\Sc^{\C}$. Moreover, extending $J$ by $0$ on $\Rc\oplus{\Rc}^{*}$, we obtain a real $2$-form $\omega$ on $M$ by setting $\omega(X,Y)=g(JX,Y)$.\\

It directly follows from formula (\ref{e12}) that 
\newpage
\begin{Le} The torsion of the canonical Witt connection $\nabla^{\C}$ on $T^{\C}M$ associated to the complex Witt structure $T^{\C}M=(\Rc^{\C}\oplus{(\Rc^{*})}^{\C})\op(T^{1,0}M\oplus T^{0,1}M)$ is given by   
\begin{eqnarray}\label{e30}
T^{\C}(\Xh,\Yh)&=&-[\Xh,\Yh]^{0,1}+d\sig{\C}{}(\Xh,\Yh)\n{}{}+d\sig{\C}{*}(\Xh,\Yh)\n{}{*}\nonumber\\
T^{\C}(\Xh,\Yah)&=&d\sig{\C}{}(\Xh,\Yah)\n{}{}+d\sig{\C}{*}(\Xh,\Yah)\n{}{*}+\Tc^{0,1}(\Xh,\Yah)-\Tc^{1,0}(\Yah,\Xh)\nonumber\\
T^{\C}(\n{}{},\Xh)&=&-[\n{}{},\Xh]^{0,1}+d\sig{\C}{*}(\n{}{},\Xh)\n{}{*}-\frac{1}{2}(\Lc_{\Xh}g^{\C})(\n{}{},\n{}{*})\n{}{}+\Tc^{1,0}(\n{}{},\Xh)\nonumber\\
T^{\C}(\n{}{*},\Xh)&=&-[\n{}{*},\Xh]^{0,1}+d\sig{\C}{}(\n{}{*},\Xh)\n{}{}-\frac{1}{2}(\Lc_{\Xh}g^{\C})(\n{}{},\n{}{*})\n{}{*}+\Tc^{1,0}(\n{}{*},\Xh)\nonumber\\
T^{\C}(\n{}{},\n{}{*})&=&-[\n{}{},\n{}{*}]_{\Sc^{\C}},
\end{eqnarray}
with 
\begin{eqnarray*}
g^{\C}(\Tc^{0,1}(\Xh,\Yah),\Zh)&=&\frac{1}{2}\Bigl((\Lc_{\Xh}g^{\C})(\Yah,\Zh)+d{(\Xh)}^{\flat}(\Yah,\Zh)\Bigr)\\
&=&\frac{\sqrt{-1}}{2}(d\omega^{\C})(\Xh,\Yah,\Zh)=\frac{1}{2}(J^{*}d\omega^{\C})(\Xh,\Yah,\Zh)\\
g^{\C}(\Tc^{1,0}(\Yah,\Xh),\Zah)&=&\frac{1}{2}\Bigl((\Lc_{\Yah}g^{\C})(\Xh,\Zah)+d{(\Yah)}^{\flat}(\Xh,\Zah)\Bigr)\\
&=&\frac{\sqrt{-1}}{2}(d\omega^{\C})(\Xh,\Yah,\Zah)=-\frac{1}{2}(J^{*}d\omega^{\C})(\Xh,\Yah,\Zah)\\ 
g^{\C}(\Tc^{1,0}(\n{}{(*)},\Xh),\Yah)&=&\frac{1}{2}(\Lc_{\n{}{(*)}}g^{\C})(\Xh,\Yah)=-\frac{\sqrt{-1}}{2}(\Lc_{\n{}{(*)}}\omega^{\C})(\Xh,\Yah),\\
-g^{\C}([\n{}{(*)},\Xh]^{0,1},\Yh)&=&-\frac{1}{2}g^{\C}((J\Lc_{\n{}{(*)}}J)(\Xh),\Yh)\\
&=&\frac{1}{2}(\Lc_{\n{}{(*)}}g^{\C})(\Xh,\Yh)+\frac{\sqrt{-1}}{2}(\Lc_{\n{}{(*)}}\omega^{\C})(\Xh,\Yh),
\end{eqnarray*}
where $\sig{\C}{},\sig{\C}{*}$ and $\omega^{\C}$ are the complex extensions of $\sig{}{},\sig{}{*}$ and $\omega$ and $\Xh,\Yh\in T^{1,0}M$, $\Yah\in T^{0,1}M$.\\
\end{Le}

\begin{Prop} 
Let $(M^{2m+2},g,\Nc)$ be an almost Robinson manifold. There exists a connection $\nabla$ on $TM$ (called the Lichnerowicz connection) such that 
\begin{enumerate}
\item $\nabla$ is a Witt connection for the Witt structure $TM=(\Rc\oplus{\Rc}^{*})\op\Sc$.
\item $\nabla J=0$ (hence $\nabla\omega=0$).
\item The torsion of $\nabla$ is given by   
\begin{eqnarray}\label{e31}
T(X,Y)&=&-\frac{1}{4}N_J(X,Y)-\frac{1}{4}\Bigl(d^{c}\omega+\Mc(d^{c}\omega)\Bigr)^{\sharp}(X,Y)+d\sig{}{}(X,Y)\n{}{}+d\sig{}{*}(X,Y)\n{}{*},\quad X,Y\in\Sc\nonumber\\
T(\n{}{},X)&=&d\sig{}{*}(\n{}{},X)\n{}{*}-\frac{1}{2}(\Lc_{X}g)(\n{}{},\n{}{*})\n{}{}+\tu{\Sc}(\n{}{},X)+\frac{1}{4}\Bigl(\Lc_{\n{}{}}\omega-\Mc(\Lc_{\n{}{}}\omega))\Bigr)^{\sharp}(JX)\nonumber\\
T(\n{}{*},X)&=&d\sig{}{}(\n{}{*},X)\n{}{}-\frac{1}{2}(\Lc_{X}g)(\n{}{},\n{}{*})\n{}{*}+\tu{\Sc}(\n{}{*},X)+\frac{1}{4}\Bigl(\Lc_{\n{}{*}}\omega-\Mc(\Lc_{\n{}{*}}\omega)\Bigr)^{\sharp}(JX)\nonumber\\
T(\n{}{},\n{}{*})&=&-[\n{}{},\n{}{*}]_{\Sc},
\end{eqnarray}
with 
\begin{eqnarray*}
N_J(X,Y)=[X,Y]-[JX,JY]+J[JX,Y]+J[X,JY],\;(d^{c}\omega)(X,Y,Z)=-(d\omega)(JX,JY,JZ),
\end{eqnarray*}
and, for $\alpha\in\Omega^{p}(M)$,
\begin{eqnarray*}
g(\alpha^{\sharp}(X_1,X_2,\ldots,X_{p-1}),X_{p})&=&\alpha(X_1,X_2,\ldots,X_{p-1},X_{p})\\
\Mc(\alpha)(X_1,X_2,X_3,\ldots,X_{p})&=&\alpha(JX_1,JX_2,X_3,\ldots,X_{p}).
\end{eqnarray*}
\end{enumerate}
\end{Prop}

\noindent{Proof.} Let $\nabla^{\C}$ be the canonical Witt connection on $T^{\C}M$ associated to the complex Witt structure $T^{\C}M=(\Rc^{\C}\oplus{(\Rc^{*})}^{\C})\op(T^{1,0}M\oplus T^{0,1}M)$. Denote by $\nabla$ the restriction of $\nabla^{\C}$ to $TM$. Since $\nabla^{\C}$ preserve $g^{\C}$ and the complex Witt structure, then $\nabla$ preserves $g$, $J$ and the real Witt structure. Now denote by $T$ the torsion of $\nabla$. Decomposing any $X\in\Sc$ as $X=\Xh+\Xah$ with $\ds\Xh=\frac{1}{2}(X-iJX)\in T^{1,0}M$ and $\Xah=\overline{\Xh}\in T^{0,1}M$, then we have, for any $X,Y\in\Sc$ 
\begin{eqnarray*}
T(X,Y)&=&T^{\C}(\Xh,\Yh)+T^{\C}(\Xh,\Yah)+T^{\C}(\Xah,\Yh)+T^{\C}(\Xah,\Yah)\\
T(\n{}{},X)&=&T^{\C}(\n{}{},\Xh)+T^{\C}(\n{}{},\Xah),\quad T(\n{}{*},X)=T^{\C}(\n{}{*},\Xh)+T^{\C}(\n{}{*},\Xah).
\end{eqnarray*}
Using formula (\ref{e30}), we obtain 
\begin{eqnarray*}
T(X,Y)&=&-[\Xh,\Yh]^{0,1}-[\Xah,\Yah]^{1,0}+\frac{1}{2}\Bigl({((J^{*}d\omega^{\C})(\Xh,\Yah))}^{\sharp}+{((J^{*}d\omega^{\C})(\Xah,\Yh))}^{\sharp}\Bigr)\\
&+&d\sig{}{}(X,Y)\n{}{}+d\sig{}{*}(X,Y)\n{}{*}.
\end{eqnarray*}
Now, we have 
$$-[\Xh,\Yh]^{0,1}-[\Xah,\Yah]^{1,0}=-\frac{1}{4}([X,Y]-[JX,JY]+J[JX,Y]+J[X,JY])=-\frac{1}{4}N_{J}(X,Y),$$
and, for $Z\in\Sc$,
\begin{eqnarray*}
\frac{1}{2}\Bigl((J^{*}d\omega^{\C})(\Xh,\Yah,Z)+(J^{*}d\omega^{\C})(\Xah,\Yh,Z)\Bigr)&=&\Re e\bigl((J^{*}d\omega^{\C})(\Xh,\Yah,Z)\bigr)\\
&=&\frac{1}{4}\bigl((d\omega)(JX,JY,JZ)+(d\omega)(X,Y,JZ)\bigr)\\
&=&-\frac{1}{4}\bigl((d^{c}\omega)(X,Y,Z)+\Mc(d^{c}\omega)(X,Y,Z)\bigr).
\end{eqnarray*}
We deduce the formula for $T(X,Y)$. Now, by (\ref{e30}), we have also 
$$
T(\n{}{},X)=-[\n{}{},\Xh]^{0,1}-[\n{}{},\Xah]^{1,0}+\Tc^{1,0}(\n{}{},\Xh)+\Tc^{0,1}(\n{}{},\Xah)+d\sig{}{*}(\n{}{},X)\n{}{*}-\frac{1}{2}(\Lc_{X}g)(\n{}{},\n{}{*})\n{}{}.
$$
We obtain
\begin{eqnarray*}
g(-[\n{}{},\Xh]^{0,1}-[\n{}{},\Xah]^{1,0}&+&\Tc^{1,0}(\n{}{},\Xh)+\Tc^{0,1}(\n{}{},\Xah),Y)\\
&=&\frac{1}{2}\Bigl((\Lc_{\n{}{}}g)(X,Y)+\sqrt{-1}(\Lc_{\n{}{}}\omega^{\C})(\Xh,\Yh)-\sqrt{-1}(\Lc_{\n{}{}}\omega^{\C})(\Xah,\Yah)\Bigr)\\
&=&\frac{1}{2}(\Lc_{\n{}{}}g)(X,Y)-\Im m\bigl((\Lc_{\n{}{}}\omega^{\C})(\Xh,\Yh)\bigr)\\
&=&\frac{1}{2}(\Lc_{\n{}{}}g)(X,Y)+\frac{1}{4}(\Lc_{\n{}{}}\omega)(JX,Y)+\frac{1}{4}(\Lc_{\n{}{}}\omega)(X,JY)\\
&=&\frac{1}{2}(\Lc_{\n{}{}}g)(X,Y)+\frac{1}{4}\Bigl(\Lc_{\n{}{}}\omega-\Mc(\Lc_{\n{}{}}\omega)\Bigr)(JX,Y).
\end{eqnarray*}
We deduce the formula for $T(\n{}{},X)$. The proof is the same for $T(\n{}{*},X)$. To conclude, we have
$T(\n{}{},\n{}{*})=T^{\C}(\n{}{},\n{}{*})=-[\n{}{},\n{}{*}]_{\Sc}$. $\Box$

\begin{Rem}
1. Note that, in general, the Lichnerowicz connection does not coincide with the canonical Witt connection associated to the real Witt structure $TM=(\Rc\oplus{\Rc}^{*})\op\Sc$.\\ 
2. For an almost Hermitian CR manifold $(M,g,J)$ (cf. example 2.2), we have a complex Witt structure given by $T^{\C}M=(T^{1,0}M\oplus T^{0,1}M)\op\epsilon(span_{\C}\{\xi\})$ (with $\xi$ nowwhere vanishing vector field on $M$ when $\epsilon=1$). In this case, we can show that the connection on $TM$ induced by the canonical Witt connection on $T^{\C}M$ is, for $\epsilon=0$, the Lichnerowicz connection (so called first canonical connection in \cite{G}) and, for $\epsilon=1$, the generalized Tanaka-Webter connection defined in \cite{NL}.
\end{Rem}

\begin{Def}
Let $(M^{2m+2},g,\Nc)$ be an almost Robinson manifold.
An almost Robinson structure will be called almost Kahler-Robinson structure if $d\omega=0$.\\
An almost Robinson structure will be called shear-free almost Robinson structure if (cf.
\cite{AGS}) 
$$\Lc_{\n{}{}}g^{\C}=\rho g^{\C}+\sig{\C}{*}\odot\alpha,$$
with $\rho\in\Cc^{\infty}(M)$ and $\alpha\in\Omega_{\C}^{1}(M)$.\\
A shear-free almost Robinson structure such that the flow 
of $\n{}{}$ acts freely on $M$ together with $\tilde{M}=M/\Fc$ manifold and $\pi:M\to\tilde{M}$ principal $G$-bundle ($G=\R$ or $S^1$) is called a regular shear-free almost Robinson structure. 
\end{Def}

An interesting subclass of regular shear-free Kahler-Robinson manifolds is given by the Fefferman spaces which are some $S^1$-principal bundles over strictly pseudoconvex $CR$ manifolds. We briefly recall the construction of these spaces given in \cite{B},\cite{L}. Let $\tilde{M}$ be a $(2m+1)$-dimensional strictly pseudoconvex $CR$ manifold. In the following, $(\tilde{\theta},\tilde{H},\tilde{\xi},\tilde{J},L_{\tilde{\theta}},\tilde{g}^{\Wc})$ respectively denote the pseudo-Hermitian structure on $\tilde{M}$, the contact hyperplan on $\tilde{M}$ (i.e. $\tilde{H}=\Ker\,\tilde{\theta}$), the Reeb field on $\tilde{M}$ (i.e. the unique vector field $\tilde{\xi}$ on $\tilde{M}$ such that $\tilde{\theta}(\tilde{\xi})=1$ and $d\tilde{\theta}(\tilde{\xi},.)=0$), the $CR$-structure on $\tilde{H}$ (i.e. the complex stucture $\tilde{J}$ on $\tilde{H}$), the Levi form (i.e. the real symmetric tensor given $L_{\tilde{\theta}}(.,.)=d\tilde{\theta}(.,\tilde{J}.)$) and the Webster metric (i.e. the riemannian metric $\tilde{g}^{\Wc}=L_{\tilde{\theta}}+\tilde{\theta}\odot\tilde{\theta}$). We assume $\tilde{M}$ endowed with its Tanaka-Webster connection $\tilde{\nabla}$ and we denote by $\tilde{R}^{\Wc}$ the curvature tensor of $\tilde{\nabla}$. The pseudo-Hermitian Ricci form $\tilde{\rho}^{\Wc}$ and the pseudo-Hermitian scalar curvature $\tilde{s}^{\Wc}$ are respectively defined by 
$$
\tilde{\rho}^{\Wc}(\tilde{X},\tilde{Y})=\sum_{i=1}^{m}\tilde{R}^{\Wc}(\tilde{X},\tilde{Y},\tilde{e}_i,\tilde{J}\tilde{e}_i) \quad\mathrm{and}\quad
\tilde{s}^{\Wc}=-\sum_{i=1}^{m}\tilde{\rho}^{\Wc}(\tilde{e}_i,\tilde{J}\tilde{e}_i),
$$
with $(\tilde{e}_1,\tilde{J}\tilde{e}_1,\ldots,\tilde{e}_m,\tilde{J}\tilde{e}_m)$ local orthonormal basis of $\tilde{H}$.\\
Let $\tilde{K}:=\wedge^{m+1,0}(\tilde{M})=\{\alpha\in\Gamma(\wedge^{m+1}T\tilde{M}^{\C}),\; i(\tilde{X}^{0,1})\alpha=0\}$ be the canonical bundle of $\tilde{M}$, $\tilde{K}^{*}=\tilde{K}\setminus\{0\}$ and $M=\tilde{K}^{*}/\R^{+}$, then $M$ is a principal $S^1$-bundle over $\tilde{M}$ with canonical projection denoted by $\pi$. Let $A^{\tilde{\nabla}}\in{\Omega}^{1}(M;i\R)$ be the connection form on $M$ associated to the Tanaka-Webster connection $\tilde{\nabla}$ on $\tilde{M}$ and let $\sig{}{*},\sig{}{}\in{\Omega}^{1}(M;\R)$ the $1$-forms on $M$ respectively given by  
$$\sig{}{*}=\pi^{*}\tilde{\theta}\quad\mathrm{and}\quad\sig{}{}=\sqrt{-1}A^{\tilde{\nabla}}+\frac{1}{2(m+1)}(\tilde{s}^{\Wc}\circ\pi)\sig{}{*}.$$ 
Note that 
\begin{equation}\label{e32}
d\sig{}{*}=\pi^{*}d\tilde{\theta}\quad\mathrm{and}\quad
d\sig{}{}=\pi^{*}\tilde{\rho}^{\Wc}+\frac{1}{2(m+1)}\Bigl(\pi^{*}d\tilde{s}^{\Wc}\wedge\sig{}{*}+(\tilde{s}^{\Wc}\circ\pi)d\sig{}{*}\Bigr).
\end{equation}
Let $\n{}{}$ (respectively $\n{}{*}$) be the fondamental vector field on $M$ satisfying $\sig{}{}(\n{}{})=1$ (respectively the vector field on $M$ such that $\pi_{*}(\n{}{*})=\tilde{\xi}$). If $\Vc_{M}=\Ker\,d\pi=span(\n{}{})$ and 
$\Hc_{M}=\Ker\,\sig{}{}=\pi_{*}^{-1}(T\tilde{M})$, then we have  
$$TM=\Vc_{M}\oplus\Hc_{M}=span(\n{}{})\oplus span(\n{}{*})\oplus\pi_{*}^{-1}(\tilde{H}).$$ 
The Fefferman metric $g^{F}$ on $M$ is  given by 
$$g^{F}=\pi^{*}L_{\tilde{\theta}}+\sig{}{*}\odot\sig{}{}.$$ 
We call $(M^{2m+2},g^{F})$ a Fefferman space. Let $\Sc=\pi_{*}^{-1}(\tilde{H})$ and $J=\pi_{*}^{-1}\circ\tilde{J}\circ\pi_{*}$, then $(g^{F}_{\Sc},J)$ is a Hermitian structure on $\Sc$. Now, $\Nc=\Vc^{\C}_{M}\oplus T^{1,0}M$ is a Robinson structure such that $\Lc_{\n{}{}}\sig{}{}=\Lc_{\n{}{}}\sig{}{*}=\Lc_{\n{}{}}g^{F}=0$ and $\omega^{F}=d\sig{}{*}$, also $(M^{2m+2},g^{F},\Nc)$ is an example of regular shear-free Kahler-Robinson manifold.\\

\begin{Prop} Let $(M^{2m+2},g^{F},\Nc)$ be a Fefferman space endowed with its Lichnerowicz connection $\nabla$, then 
\begin{enumerate}
\item A curve $c$ on $M$ is a $\nabla$-lightlike geodesic if and only if $(\Lc_{\dot{c}}\sig{}{})(\dot{c})=(\Lc_{\dot{c}}\sig{}{*})(\dot{c})=0$.
\item A curve $c\in\Cc_{\Sc}$ is a normal sub-riemannian geodesic if and only if
\begin{eqnarray}\label{e33}
\nabla_{\dot{c}}\dot{c}=k_1{\bigl(i(\dot{c})\pi^{*}\tilde{\rho}^{\Wc}\bigr)}_{\Sc}^{\sharp}&+&\Bigl(k_1\Bigl(\frac{\tilde{s}^{\Wc}\circ\pi}{2(m+1)}+\int_{0}^{t}(\pi^{*}\tilde{\rho}^{\Wc})(\n{}{*},\dot{c})\,ds-\frac{1}{2(m+1)}\int_{0}^{t}(\pi^{*}d\tilde{s}^{\Wc})(\dot{c})\,ds\Bigr)\nonumber\\
&&+\frac{1}{2}\int_{0}^{t}(\pi^{*}\Lc_{\tilde{\xi}}g^{\Wc})(\dot{c},\dot{c})\,ds+k_2\Bigr)J\dot{c}.
\end{eqnarray}
In particular, if $(M,g^{F})$ is a Fefferman space over a pseudo-Einstein Sasaki manifold, then the last equation becomes 
$$\nabla_{\dot{c}}\dot{c}=\Bigl(-k_1\frac{m+2}{2m(m+1)}(\tilde{s}^{\Wc}\circ\pi)+k_2\Bigr)J\dot{c}.$$
\end{enumerate}
\end{Prop}

\noindent{Proof.} If $(M,g^{F},\Nc)$ is a Fefferman space, then we have the relations  
\begin{equation}\label{e34}
N_{J}=0,\quad \omega^{F}=d\sig{}{*},\quad \Lc_{\n{}{}}\sig{}{}=\Lc_{\n{}{}}\sig{}{*}=\Lc_{\n{}{}}g^{F}=\Lc_{\n{}{}}\omega^{F}=0,\quad \Lc_{\n{}{*}}\sig{}{*}=\Lc_{\n{}{*}}\omega^{F}=0.
\end{equation}
Substituting these relations in (\ref{e31}), we deduce that the Lichnerowicz connection coincides with the canonical Witt connection associated to the Witt structure $TM=(\Rc\oplus{\Rc}^{*})\op\Sc$. Hence, the equations for $\nabla$-lightlike geodesics follow from (\ref{e14}) together with $\Lc_{\n{}{}}\sig{}{}=\Lc_{\n{}{*}}\sig{}{*}=0$. 
Now, let $c\in\Cc_{\Sc}$ be a normal sub-riemannian geodesic, formula (\ref{e22}) together with relations (\ref{e34}) gives 
\begin{eqnarray}\label{e35}
\nabla_{\dot{c}}\dot{c}&=&\lambda_{1}^{c}{\bigl(i(\dot{c})d\sig{}{}\bigr)}_{H}^{\sharp}+\lambda_{2}^{c}{\bigl(i(\dot{c})d\sig{}{*}\bigr)}_{H}^{\sharp}\nonumber\\
\begin{pmatrix}\dot{\lambda_{1}^{c}}\cr\dot{\lambda_{2}^{c}}\end{pmatrix}&=&\begin{pmatrix}0\cr\lambda_{1}^{c}(\Lc_{\n{}{*}}\sig{}{})(\dot{c})+\frac{1}{2}(\Lc_{\n{}{*}}g^{F})(\dot{c},\dot{c})\end{pmatrix}.
\end{eqnarray}
Using (\ref{e32}), we have, for $X_{\Sc},Y_{\Sc}\in\Sc$ 
\begin{eqnarray*}
d\sig{}{}(X_{\Sc},Y_{\Sc})&=&(\pi^{*}\tilde{\rho}^{\Wc})(X_{\Sc},Y_{\Sc})+\frac{1}{2(m+1)}(\tilde{s}^{\Wc}\circ\pi)\omega^{F}(X_{\Sc},Y_{\Sc}),\\
(\Lc_{\n{}{*}}\sig{}{})(X_{\Sc})&=&d\sig{}{}(\n{}{*},X_{\Sc})=(\pi^{*}\tilde{\rho}^{\Wc})(\n{}{*},X_{\Sc})-\frac{1}{2(m+1)}(\pi^{*}d\tilde{s}^{\Wc})(X_{\Sc}),\\
(\Lc_{\n{}{*}}g^{F})(X_{\Sc},Y_{\Sc})&=&(\pi^{*}\Lc_{\tilde{\xi}}g^{\Wc})(X_{\Sc},Y_{\Sc}).
\end{eqnarray*}
Hence (\ref{e35}) becomes 
\begin{eqnarray*}
\nabla_{\dot{c}}\dot{c}&=&\lambda_{1}^{c}\bigl({i(\dot{c})\pi^{*}\tilde{\rho}^{\Wc}\bigr)}_{\Sc}^{\sharp}+
\Bigl(\frac{\lambda_{1}^{c}(\tilde{s}^{\Wc}\circ\pi)}{2(m+1)}+\lambda_{2}^{c}\Bigr)J\dot{c}\\
\begin{pmatrix}\dot{\lambda_{1}^{c}}\cr\dot{\lambda_{2}^{c}}\end{pmatrix}&=&\begin{pmatrix}0\cr\ds\lambda_{1}^{c}(\pi^{*}\tilde{\rho}^{\Wc})(\n{}{*},\dot{c})-\frac{\lambda_{1}^{c}}{2(m+1)}(\pi^{*}d\tilde{s}^{\Wc})(\dot{c})+\frac{1}{2}(\pi^{*}\Lc_{\tilde{\xi}}g^{\Wc})(\dot{c},\dot{c})\end{pmatrix}.
\end{eqnarray*}
Hence $\lambda_{1}^{c}$ is constant equal to $k_1$ and by integrating we obtain :
$$\lambda_{2}^{c}=k_1\int_{0}^{t}(\pi^{*}\tilde{\rho}^{\Wc})(\n{}{*},\dot{c})-\frac{1}{2(m+1)}(\pi^{*}d\tilde{s}^{\Wc})(\dot{c})\,ds+\frac{1}{2}\int_{0}^{t}(\pi^{*}\Lc_{\tilde{\xi}}g^{\Wc})(\dot{c},\dot{c})\,ds+k_2\quad (k_2\in\R).$$
Formula (\ref{e33}) directly follows. Now if $\tilde{M}$ is a pseudo-Einstein Sasaki manifold then we have 
$(\Lc_{\tilde{\xi}}g^{\Wc})=0$ and $\ds\tilde{\rho}^{\Wc}=-\frac{\tilde{s}^{\Wc}}{m}d\tilde{\theta}$ with $\tilde{s}^{\Wc}$ constant. Hence $\ds(\pi^{*}\tilde{\rho}^{\Wc})=-\frac{\tilde{s}^{\Wc}\circ\pi}{m}\omega^{F}$ and we obtain the result in this case. $\Box$\\

Let $(M,g,\Nc)$ be an almost Robinson manifold endowed with its complex Witt structure $T^{\C}M=(\Rc^{\C}\oplus{(\Rc^{*})}^{\C})\op(T^{1,0}M\oplus T^{0,1}M)$. By setting $V_{+}=\Rc^{\C}\oplus{(\Rc^{*})}^{\C}$ and $V_{-}=T^{1,0}M\oplus T^{0,1}M=\Sc^{\C}$, we define the isometric involution $\delta_{x}:T^{\C}_{x}M\to T^{\C}_{x}M$  by $\ds\delta_{x}=\pi^{+}_{x}-\pi^{-}_{x}$ where $\pi^{\pm}_{x}$ are the projections onto $V_{\pm,x}$.

\begin{Def} Let $(M^{2m+2},g,\Nc)$ be an almost Robinson manifold endowed with its complex Witt structure $T^{\C}M=(\Rc^{\C}\oplus{(\Rc^{*})}^{\C})\op(T^{1,0}M\oplus T^{0,1}M)$ and its canonical Witt connection $\nabla^{\C}$, then $(M,g,\Nc)$ is called an (locally) almost Robinson Witt symmetric space if, for any point $x\in M$, $\psi^{\C}_{x}=exp^{\nabla^{\C}}_{x}\circ\delta_{x}\circ {(exp^{\nabla^{\C}}_{x})}^{-1}$ is a (local) Witt isometry. 
\end{Def}

\begin{Prop} Let $(M^{2m+2},g,\Nc)$ be a locally almost Robinson Witt symmetric space endowed with its Lichnerowicz connection $\nabla$, then we have 
\begin{enumerate}
\item There is a codimension $2m$ foliation $\Fc$ of $M$ such that $T(\Fc)=\Rc\oplus{\Rc}^{*}$, 
\item The almost Hermitian structure $(g_{\Sc},J)$ is integrable and ${(d\omega^{\C})}_{\Sc^{\C}}^{(2,1)+(1,2)}=0$,
\item ${(\Lc_{\n{}{}}\sig{}{})}_{\Sc}=-{(\Lc_{\n{}{*}}\sig{}{*})}_{\Sc}$, ${(\Lc_{\n{}{}}\sig{}{*})}_{\Sc}={(\Lc_{\n{}{*}}\sig{}{})}_{\Sc}=0$,\\
\\
${(\nab{\Sc}d\sig{}{})}_{\Sc}=-{(\Lc_{\n{}{}}\sig{}{})}_{\Sc}\otimes {(d\sig{}{})}_{\Sc}$, ${(\nab{\Sc}d\sig{}{*})}_{\Sc}=-{(\Lc_{\n{}{*}}\sig{}{*})}_{\Sc}\otimes {(d\sig{}{*})}_{\Sc}$,\\
\\
${(\nab{\Sc}\Lc_{\n{}{}}g)}_{\Sc}={(\nab{\Sc}\Lc_{\n{}{*}}g)}_{\Sc}={(\nab{\Sc}\Lc_{\n{}{}}\omega)}_{\Sc}={(\nab{\Sc}\Lc_{\n{}{*}}\omega)}_{\Sc}=0$,
\item $\nab{\Sc}R=0$, ${(i_{\n{}{}}R)}_{\Sc}={(i_{\n{}{*}}R)}_{\Sc}=0$.
\end{enumerate}
\end{Prop} 

\noindent{Proof.} If $M$ is a locally almost Robinson Witt symmetric space, then, by proposition 5.1, we have $[V_{+},V_{+}]\subset V_{+}$. Hence $V_{+}$ is an integrable $2$-dimensional distribution and we have a foliation $\Fc$ of $M$ with $T(\Fc)=\Rc\oplus{\Rc}^{*}$. Now, we have also that $T^{\C}$ must satisfies $T^{\C}_{/V_{-}\times V_{-}}\subset V_{+}$ and $T^{\C}_{/V_{+}\times V_{-}}\subset V_{-}$. These conditions applied to (\ref{e30}) yields to 
\begin{equation}\label{e36}
[\Xh,\Yh]^{0,1}=[\Xah,\Yah]^{1,0}=0,\quad{((J^{*}d\omega^{\C})(\Xh,\Yah))}^{\sharp}=0,
\end{equation}
\begin{eqnarray}\label{e37}
d\sig{\C}{*}(\n{}{},\Xh)&=&d\sig{\C}{*}(\n{}{},\Xah)=0,\quad d\sig{\C}{}(\n{}{*},\Xh)=d\sig{\C}{}(\n{}{*},\Xah)=0,\nonumber\\
(\Lc_{\Xh}g^{\C})(\n{}{},\n{}{*})&=&(\Lc_{\Xah}g^{\C})(\n{}{},\n{}{*})=0.
\end{eqnarray}
It directly follows from (\ref{e36}) that $J$ is integrable (i.e. $N_J=0$) and ${(d\omega^{\C})}_{\Sc^{\C}}^{(2,1)+(1,2)}=0$ (since ${((J^{*}d\omega^{\C})(\Xh,\Yah))}^{\sharp}=0$). Now, for $X=\Xh+\Xah\in\Sc$, we have by (\ref{e37}) 
$$d\sig{}{*}(\n{}{},X)=(\Lc_{\n{}{}}\sig{}{*})(X)=0,\; d\sig{}{}(\n{}{*},X)=(\Lc_{\n{}{*}}\sig{}{})(X)=0,\;(\Lc_{X}g)(\n{}{},\n{}{*})=(\Lc_{\n{}{}}\sig{}{})(X)+(\Lc_{\n{}{*}}\sig{}{*})(X)=0.
$$
By proposition 5.1, we have also $\nabla^{\C}_{\Sc^{\C}}T^{\C}=0$. We first deduce that 
$${\bigl(\nabla^{\C}_{\Sc^{\C}}\Lc_{\n{}{}}g^{\C}\bigr)}_{\Sc^{\C}}={\bigl(\nabla^{\C}_{\Sc^{\C}}\Lc_{\n{}{*}}g^{\C}\bigr)}_{\Sc^{\C}}={\bigl(\nabla^{\C}_{\Sc^{\C}}\Lc_{\n{}{}}\omega^{\C}\bigr)}_{\Sc^{\C}}={\bigl(\nabla^{\C}_{\Sc^{\C}}\Lc_{\n{}{*}}\omega^{\C}\bigr)}_{\Sc^{\C}}=0.$$
Now,   
\begin{eqnarray*}
{\bigl(\nabla^{\C}_{\Sc^{\C}}(d\sig{\C}{}\otimes\n{}{})\bigr)}_{\Sc^{\C}}&=&{\bigl(\nabla^{\C}_{\Sc^{\C}}d\sig{\C}{}\bigr)}_{\Sc^{\C}}\otimes\n{}{}+{(d\sig{\C}{})}_{\Sc^{\C}}\otimes\nabla^{\C}_{\Sc^{\C}}\n{}{}\\
&=&\Bigl({\bigl(\nabla^{\C}_{\Sc^{\C}}d\sig{\C}{}\bigr)}_{\Sc^{\C}}-\bigl(\nabla^{\C}_{\Sc^{\C}}\sig{\C}{}\bigr)(\n{}{}){(d\sig{\C}{})}_{\Sc^{\C}}\Bigr)\otimes\n{}{}\\
&=&\Bigl({\bigl(\nabla^{\C}_{\Sc^{\C}}d\sig{\C}{}\bigr)}_{\Sc^{\C}}-\frac{1}{2}{\bigl(\Lc_{\n{}{*}}\sig{\C}{*}-\Lc_{\n{}{}}\sig{\C}{}\bigr)}_{\Sc^{\C}}\otimes{(d\sig{\C}{})}_{\Sc^{\C}}\Bigr)\otimes\n{}{}\\
&=&\Bigl({\bigl(\nabla^{\C}_{\Sc^{\C}}d\sig{\C}{}\bigr)}_{\Sc^{\C}}+{\bigl(\Lc_{\n{}{}}\sig{\C}{}\bigr)}_{\Sc^{\C}}\otimes{(d\sig{\C}{})}_{\Sc^{\C}}\Bigr)\otimes\n{}{},
\end{eqnarray*}
and 
$${\bigl(\nabla^{\C}_{\Sc^{\C}}(d\sig{\C}{*}\otimes\n{}{*})\bigr)}_{\Sc^{\C}}
=\Bigl({\bigl(\nabla^{\C}_{\Sc^{\C}}d\sig{\C}{*}\bigr)}_{\Sc^{\C}}+{\bigl(\Lc_{\n{}{*}}\sig{\C}{*}\bigr)}_{\Sc^{\C}}\otimes{(d\sig{\C}{*})}_{\Sc^{\C}}\Bigr)\otimes\n{}{*}.$$
The assumption $\nabla^{\C}_{\Sc^{\C}}T^{\C}=0$ implies that 
$$
{\bigl(\nabla^{\C}_{\Sc^{\C}}d\sig{\C}{}\bigr)}_{\Sc^{\C}}=-{\bigl(\Lc_{\n{}{}}\sig{\C}{}\bigr)}_{\Sc^{\C}}\otimes{(d\sig{\C}{})}_{\Sc^{\C}}\quad{\mathrm and}\quad{\bigl(\nabla^{\C}_{\Sc^{\C}}d\sig{\C}{*}\bigr)}_{\Sc^{\C}}=-{\bigl(\Lc_{\n{}{*}}\sig{\C}{*}\bigr)}_{\Sc^{\C}}\otimes{(d\sig{\C}{*})}_{\Sc^{\C}}.
$$
Parallelism for the real tensors directly follows.\\
To conclude, since $\nabla^{\C}_{\Sc^{\C}}R^{\C}=0$ and $R^{\C}_{/V_{+}\times V_{-}}=0$, then $\nab{\Sc}R=0$ and ${(i_{\n{}{}}R)}_{\Sc}={(i_{\n{}{*}}R)}_{\Sc}=0$. $\Box$\\

For a Fefferman Witt symmetric space, the previous proposition becomes 

\begin{Cor} Let $(M^{2m+2},g^{F})$ be a locally Fefferman Witt symmetric space endowed with its Lichnerowicz connection $\nabla$, then we have  
\begin{enumerate}
\item There is a codimension $2m$ foliation $\Fc$ of $M$ such that $T(\Fc)=\Rc\oplus{\Rc}^{*}$, 
\item $\Lc_{\n{}{}}\sig{}{}=\Lc_{\n{}{*}}\sig{}{}=\Lc_{\n{}{}}\sig{}{*}=\Lc_{\n{}{*}}\sig{}{*}=0,\quad \nab{}\sig{}{}=\nab{}\sig{}{*}=\nab{}d\sig{}{*}=0,\quad {(\nab{\Sc}d\sig{}{})}_{\Sc}={(\nab{\Sc}\Lc_{\n{}{*}}g^{F})}_{\Sc}=0,$
\item $\nab{\Sc}R=0$, ${(i_{\n{}{}}R)}_{\Sc}={(i_{\n{}{*}}R)}_{\Sc}=0$.
\end{enumerate}
\end{Cor}

\end{document}